\documentclass[preprint,12pt]{elsarticle}

\usepackage{color}
\usepackage{graphicx}
\usepackage{latexsym,amsbsy,amssymb,amsmath}
\usepackage{graphicx}
\usepackage[justification=centering]{caption}
\usepackage{subcaption}
\usepackage{geometry}
\usepackage{enumerate}
\usepackage[Sonny]{fncychap}
\usepackage{color}
\usepackage{geometry}
\usepackage{graphicx}
\usepackage{bibtopic}
\usepackage[]{natbib}
\journal{ }

\begin{document}
\begin{frontmatter}

\title{The Penalty Cell-Centered Finite Element Scheme For Stokes Problem On General Meshes}        
\author[hai]{Ong Thanh Hai}
\author[phuong]{T.T.P. Hoang}
\author[hai] {H. Nguyen Xuan}   
\address[hai]{Faculty of Mathematics and Computer Science, University of Science, VNU HCMC, 227 Nguyen Van Cu Street, District 5, Ho Chi Minh City, Vietnam.}  
\address[phuong]{Faculty of  Mathematics-Informatics, Ho Chi Minh City University of Pedagogy, 280 An Duong Vuong Street, District 5, Ho Chi Minh City, Vietnam.}
\date{25 Mar 2014}         

\begin{abstract}
The paper is devoted to the penalty cell-centered finite element scheme (pFECC) on general meshes for the stationary Stokes problems with an incompressible variable viscosity and Dirichlet boundary conditions. In the objectives of this work, we show the rigorous mathematical analysis including the existence, the uniqueness of a discrete solution of the problem, the symmetric and the positive definite stiffness matrix, convergence of the pFECC scheme.
\end{abstract}
\begin{keyword}
The cell-centered finite element scheme; Penalty method; General grids; Volumetric locking; The stationary Stokes equations for an incompressible variable viscous fluid.
\end{keyword}
\end{frontmatter}
\section{Introduction}
\label{sec1}
Let $\Omega$ be a open, bounded domain of $\mathbb{R}^2$ with the boundary $\partial \Omega$. We consider the stationary Stokes problem for an incompressible variable viscosity in $\Omega$: find an approximation weak solution of $\bold u = (u^{(1)},u^{(2)}) \in (H^1_0(\Omega))^2$ and $p \in L^2(\Omega)$,  to the following problem
\begin{eqnarray}
\label{ptr1}
- \text{div} (2 \mu(\bold x) \boldsymbol{\mathcal E}(\bold u)) + \nabla p &=& \bold f, \quad \text{in}~\Omega, \nonumber \\
\text{div} (\bold u) &=& 0, \quad \text{in}~\Omega,\\
\bold u &=& 0, \quad \text{on}~\partial\Omega \nonumber,
\end{eqnarray}
where the velocity $\bold u: \Omega \to \mathbb{R}^2$, has the two components $u^{(1)},~u^{(2)}$, the pressure $p$ defined over $\Omega$, the strain tensor related to the displacement is defined by $\boldsymbol{\mathcal E}(\bold u) = \frac{1}{2}(\bold {\nabla u} + \bold{\nabla}^T \bold{u} )$ and $\bold f$ are the body forces per unit mass.\\
In the physical model (\ref{ptr1}), the variable viscosity (the measure of a fluid's ability to resist gradual deformation by shear or tensile stresses) of non-Newtonian fluids is dependent on shear rate or shear rate history. This interesting physical model is appeared in many commonly found substances such as ketchup, custard, toothpaste, starch suspensions, paint, blood, and shampoo...In this paper, we also assume that the viscosity $\mu: \Omega \to \mathbb R$ in (\ref{ptr1}) is piecewise Lipschitz-continuous on the domain $\overline\Omega$ and there exists $\lambda, \underline\lambda, \overline \lambda$ such that 
\begin{equation}
\label{ptr2}
\underline \lambda \le \mu(\bold x) \le \overline  \lambda, \quad \text{for a.e $\bold x \in \Omega$}
\end{equation}
and
\begin{equation}
\label{ptr3}
|\mu(\bold x) - \mu( \overline{\bold x})| \le \lambda |\bold x -   \overline{\bold x}|, 
\quad \text{for all $\bold x, \overline{\bold x} \in \overline \Omega$}
\end{equation}
With the important physical role of the Stokes equations (\ref{ptr2}), there are many numerical schemes have been extensively studied: see [11, 20, 21, 22, 13, 12] and references therein. Among different schemes, finite element schemes andfinite volume schemes are frequently used for mathematical or engineering studies.
\section{Stokes problems}
\label{sec1a}
Under hypotheses $\bold f = (f^{(1)},f^{(2)}) \in (L^2(\Omega))^2$, (\ref{ptr2}) and (\ref{ptr3}), let
\begin{equation}
\label{ptr4}
\mathbb H(\Omega) = \left\{\bold v = (v^{(1)},v^{(2)}) \in (H^1_0(\Omega))^2, \text {div}(\bold v) =   \frac{{\partial {v^{(1)}}}}{{\partial {x_1}}} + \frac{{\partial {v^{(2)}}}}{{\partial {x_2}}} = 0 \right\},
\end{equation}
then the weak solution $\bold u = (u^{(1)},u^{(2)})$ of (\ref{ptr1}) (see e.g. \cite{Ab,Te}) must be satisfied 
\begin{equation}
\label{weakprt1}
\left\{ \begin{array}{l}
\bold u = ({u^{(1)}},{u^{(2)}}) \in \mathbb H(\Omega ),\\
\int\limits_\Omega  {\eta (\bold x)\nabla \bold u(\bold x):\nabla \bold v(\bold x) \text{d}\bold x}  = \int\limits_\Omega  {\bold {f(x)v(x)} \text{d}\bold x} ,\,\forall \bold v = ({v^{(1)}},{v^{(2)}}) \in \mathbb H(\Omega)
\end{array} \right.
\end{equation}
with $\bold x = (x_1,x_2)$ and $\bold{\nabla u(x) : \nabla v(x)} = \nabla u^{(1)}(\bold x).\nabla v^{(1)}(\bold x) + \nabla u^{(2)}(\bold x).\nabla v^{(2)}(\bold x)$.\\
\newline
In order to study convergence of the approximate solution, we need the regularity of the weak solution $(\bold u, p)$. Thank to Lemma 5.2.5 of \cite{Ab}, the author proved that if the viscosity $\mu$ belongs to $C^2(\Omega)$, then the solution $(\bold u, p)$ satisfy
\begin{equation}
\label{ptr5}
\bold u \in (H^2(\Omega))^2~\text{and}~p \in H^1(\Omega).
\end{equation}
\section{The penalty cell-centered finite element framework}
\label{sec3}
The cell-centered finite element scheme (FECC), which was firstly introduced by Christophe and Ong \citep{PoHa}, was applied into the diffusion problems on general meshes. To develop the idea of the scheme for Stokes problems, we combine the FECC and the stabilization inspired by the well-known penalty method \cite{BrePit, OdKiSo} in the finite element framework.
\subsection{Discretization of the domain $\Omega$}
\label{subsec1}
For a given $\Omega$ be an open bounded polygonal set of $\mathbb{R}^2$ with boundary $\partial \Omega$. In order to partition the domain $\Omega$, we use the three families $\mathcal{D}(\mathcal{M},\mathcal E,\mathcal{P})$, $\mathcal{D^{*}}(\mathcal{M^{*}},\mathcal E^*,\mathcal{P^*})$ and $\mathcal{D^{**}}(\mathcal{M}^{**},\mathcal E^{**},\mathcal{P}^{**})$ constructed in \cite{PoHa}. The first mesh $\mathcal M$ is assumed that any line which connects two mesh points of  two adjacent elements of $\mathcal M^*$ intersects with the common edge of these elements at the unique point. Without loss of generality, we can choose each dual mesh point of $\mathcal M^*$ located at a vertex of $\mathcal M$. We introduce elements of $\mathcal M^{**}$ denoted by $L^{**}$ or $K^{**}$ or $T_{K^*,\sigma}$, where let an edge $\sigma$ of $K^* \in \mathcal M^*$, a triangle $T_{K^*,\sigma}$ stays in $K^*$.
\subsubsection{Geometrical Conditions For Regular Meshes}
\label{family4}
To study the convergence of the scheme, we will need the size of the third mesh $\mathcal M^{**}$ defined by 
\begin{equation}
\label{ptr6}
h = \mathop {\sup }\limits_{K^{**} \in \mathcal M^{**}} ~\text{diam}(K^{**}),
\end{equation}
where $\text{diam}(K^*)$ indicates a diameter of the circumscribed circle of $K^{**}$.\\
And the regularity of the three meshes is required by the existence of positive real numbers $C_i$ such that
\begin{equation}
\label{ptr7}
\text{card}(\mathcal E_{K^*}) \le C_1,~\quad \text{for all}~K^*\in\mathcal M^*,
\end{equation}
where $\mathcal E^*_{K^*}$ is a set of all edges of $K^*$. This corresponds to  the condition 
\begin{displaymath}
\text{card}(\mathcal M^{**}_{K^*}) \le C_1~\quad \text{for all}~K^*\in\mathcal M^*
\end{displaymath}
in which $\mathcal M^{**}_{K^*}$ is a set of all element of $\mathcal M^{**}$ staying in $K^{**}$.
\begin{equation}
\label{ptr10}
\text{diam}^2(K^*) \le C_2 {m}(K^*), \forall K^* \in \mathcal{M^*}.
\end{equation}
\begin{equation}
\label{ptr10a}
\text{diam}^2(T) \le C_3 {m}(L^{**}), \forall L^{**} \in \mathcal{M^{**}}.
\end{equation}
Besides, we also have another condition for the third mesh $\mathcal M^{**}$:\\
\newline
\textit{Inverse assumption: There exists constant $\zeta_{\mathcal D^{**}} > 0$ such that}
\begin{equation}
\label{ptr10b}
\mathop {\max }\limits_{K^{**} \in {\mathcal M^{**}}} \frac{{{h}}}{{\text{diam}(K^{**})}} \le \zeta_{\mathcal D^{**}} \quad \textit{for all}~h > 0.
\end{equation}
\subsection{Unknowns and Discrete operators}
\label{family5}
We will express the new scheme in the weak form; to this aim, let us firstly define the sets containing the discrete unknowns, the discrete operators, the discrete gradient, and the discrete divergence:\\
\newline
For given two neighboring elements of the first mesh $\mathcal M$, we assume that the line joining their primary mesh points can be intersected their common edge. With this assumption, let $\sigma \in \mathcal E_\text{int}~\text{such that}~\mathcal{M}_\sigma=\{K,L\}$, the three points $\bold x_K$, $\bold x_L \in \mathcal{P} ~\text {and}~ \bold x_{K^*} \in \mathcal{P^*}$ can generate a triangular element $(\bold x_K, \bold x_L, \bold x_{K^*})$  of $\mathcal M^{**}$. On this triangle, we take the unknown values $\bold u_{K^*}= (u^{(1)}_{K^*},u^{(2)}_{K^*}),~\bold u_K = (u^{(1)}_{K},u^{(2)}_{K}),~\bold u_L = (u^{(1)}_{L}),u^{(2)}_{L})$ of the velocity $\bold u = (u^{(1)},u^{(2)})$ at $\bold x_{K^*},~\bold x_K,~\bold x_L$. Besides, we have a notation $\bold u^{K^*}_{\sigma}$ (a temporary unknown) seen as a value of $\bold u $ at $\bold x_{\sigma}$, where the point $\bold x_\sigma$ is an intersecting point between the line joining two mesh points $\bold x_K$, $\bold x_L$ and the internal edge $\sigma$.\\
From these values, we introduce the following discrete velocity space\\
\newline
{\bf Definition \ref{sec3}.1}: {\em Let us define the discrete function space $\mathcal H_{\mathcal D}$ as the set of all
$((\bold u_K)_{K \in \mathcal M},(\bold u_{K^*})_{K^* \in \mathcal M^*})$, $\bold u_{K} \in \mathbb R^2$ for all $K \in \mathcal M$ and $\bold u_{K^*} \in \mathbb R^2$ for all $K \in \mathcal M^*$. Moreover, the value $\bold u_{K^*}$ is equal to 0, while a mesh point $\bold x_{K^*}$ belongs to the boundary $\partial \Omega$}.\\
\newline
and the discrete pressure space $\mathcal L_{\mathcal D}$:\\
\newline
{\bf Definition \ref{sec3}.2}:  {\em The space $\mathcal{L}_{\mathcal D}$ contains all piecewise constant functions on the dual mesh $\mathcal M^*$.}
\begin{displaymath}
\mathcal L_{\mathcal D} =  \left \{q: \Omega \to \mathbb R | q(\bold x) = \sum\limits_{{K^*} \in {\mathcal M^*}} {{q_{{K^*}}}{\chi _{{K^*}}}({\bf{x}})} \right \},
\end{displaymath}
with the characteristic function  $\chi_{K^*}$, for each $K^* \in \mathcal M^*$.\\
\newline
From the definition of the two discrete spaces $\mathcal H_{\mathcal D}$ and $\mathcal L_{\mathcal D}$, we construct a discrete gradient $\nabla_{\mathcal D,\Lambda} u_h^{(i)}$ and the interpolation operator $P_{\left( {K^* ,K,L} \right)} (u_h^{(i)})$, $i = 1, 2$, on two sub-triangles of the triangle $(\bold x_{K^*}, \bold x_K, \bold x_L)$, where their definitions are taken into account the variable viscosity $\mu(\bold x)$ and an element $\bold u_h = ((\bold u_K)_{K \in \mathcal M},(\bold u_{K^*})_{K^* \in \mathcal M^*})\in \mathcal H_{\mathcal D}$, as follows:\\
The interpolation operator for each its element $u_h^{(i)}$, i = 1,2, is defined by
\begin{displaymath}
P_{\left( {K^* ,K,L} \right)} (u_h^{(i)})~:\quad \left( {\bold x_{K^* } ,\bold x_K ,\bold x_L } \right) \to \mathbb R,
\end{displaymath}
such that it is continuous, piecewise linear on $(\bold x_{K^*},\bold x_K,\bold x_{\sigma})$ and $(\bold x_{K^*}, \bold x_L, \bold x_{\sigma})$ (two sub-triangles of $(\bold x_K, \bold x_L,\bold x_{K^*})$).\\
\begin{itemize}
{\item on the sub-triangle $(\bold x_{K^*},\bold x_K,\bold x_{\sigma})$
\begin{displaymath}
P_{(K^* ,K,L)} (u_h^{(i)})(\bold x) = \left\{ \begin{array}{l}
 u^{(i)}_K \quad \bold x = \bold x_K,  \\ 
 u^{(i)}_{K^* } \quad \bold x = \bold x_{K^* },  \\ 
 u_{\sigma,K^*} ^{(i)} \quad \bold x = \bold x_\sigma.   \\ 
 \end{array} \right.
\end{displaymath}
\begin{eqnarray*}
&&\nabla _{\mathcal D,\mu } u^{(i)} = \nabla _{\mathcal D,\mu } P_{\left( {K^* ,K,L} \right)}(u^{(i)}) \\
&=& \frac{{ - P_{\left( {K^* ,K,L} \right)}(u_h^{(i)})(\bold x_\sigma  ) n_{[\bold x_{K^*} ,\bold x_K]} - P_{\left( {K^* ,K,L} \right)}(u_h^{(i)})(\bold x_K ) n^K_{[\bold x_\sigma, \bold x_{K^*}]} - P_{\left( {K^* ,K,L} \right)}(u_h^{(i)})(\bold x_{K^* } ) n_{[\bold x_\sigma, \bold x_K]}}}{{2m_{\left( {\bold x_{K^* } ,\bold x_K ,\bold x_\sigma  } \right)} }}\\
&=& \frac{{ - u^{(i)}_{\sigma,K^*}n_{[\bold x_{K^*} ,\bold x_K]}   - u^{(i)}_{K} n^K_{[\bold x_\sigma, \bold x_{K^*}]}   -  u^{(i)}_{K^* } n_{[\bold x_\sigma, \bold x_K]} }}{{2m_{\left( {\bold x_{K^* } ,\bold x_K ,\bold x_\sigma  } \right)} }},
\end{eqnarray*}
where $n^K_{[\bold x_\sigma, \bold x_{K^*}]}$ is outer normal vector to the triangle $(\bold x_{K^*},\bold x_K,\bold x_{\sigma})$. The length of vector $n^K_{[\bold x_\sigma, \bold x_{K^*}]}$ is equal to the length of segment $[\bold x_\sigma, \bold x_{K^*}]$. If $\bold x_\sigma$ belongs to boundary $\partial \Omega$ then $u^{(i)}_{\sigma,K^*}=0$. A notation $m_{(\bold x_{K^*},\bold x_K, \bold x_\sigma)}$ is the area of a triangle $(\bold x_{K^*},\bold x_K, \bold x_\sigma)$}
{\item on the sub-triangle $(\bold x_{K^*},\bold x_L,\bold x_{\sigma})$
\begin{displaymath}
P_{(K^* ,K,L)} (u_h^{(i)})(\bold x) = \left\{ \begin{array}{l}
 u^{(i)}_L \quad \bold x = \bold x_K,  \\ 
 u^{(i)}_{K^* } \quad \bold x = \bold x_{K^* },  \\ 
 u_{\sigma,K^*}^{(i)} \quad \bold x = \bold x_\sigma.  \\ 
 \end{array} \right.
\end{displaymath}
\begin{eqnarray*}
&&\nabla _{\mathcal D,\mu } u_h^{(i)} = \nabla _{\mathcal D,\mu } P_{\left( {K^* ,K,L} \right)}(u_h^{(i)}) \\
&=& \frac{{ - P_{\left( {K^* ,K,L} \right)}(u_h^{(i)})(\bold x_\sigma  ) n_{[\bold x_{K^*} ,\bold x_L]}  - P_{\left( {K^* ,K,L} \right)}(u_h^{(i)})(\bold x_L ) n^L_{[\bold x_\sigma, \bold x_{K^*}]}  -  P_{\left( {K^* ,K,L} \right)}(u_h^{(i)})(\bold x_{K^* } ) n_{[\bold x_\sigma, \bold x_L]}}}{{2m_{\left( {\bold x_{K^* } ,\bold x_L ,\bold x_\sigma  } \right)} }}\\
& = & \frac{{ - u^{(i)}_{\sigma,K^* } n_{[\bold x_{K^*} ,\bold x_L]}   - u^{(i)}_L n^L_{[\bold x_\sigma, \bold x_{K^*}]}   - u^{(i)}_{K^* } n_{[\bold x_\sigma, \bold x_L]}}}{{2m_{\left( {\bold x_{K^* } ,\bold x_L ,\bold x_\sigma  } \right)} }},
\end{eqnarray*}
where $n^L_{[\bold x_\sigma, \bold x_{K^*}]}$ is outer normal vector to the triangle $(\bold x_{K^*},\bold x_L,\bold x_{\sigma})$. The length of vector $n^L_{[\bold x_\sigma, \bold x_{K^*}]}$ equal to the length of segment $[\bold x_\sigma, \bold x_{K^*}]$.  A notation $m_{(\bold x_{K^*},\bold x_K, \bold x_\sigma)}$ is the area of a triangle $(\bold x_{K^*},\bold x_K, \bold x_\sigma)$.}
\end{itemize}
These definitions depend on $u^{(i)}_{\sigma,K^*}$, but this temporary unknown can be fixed by imposing the Local Conservativity of the Fluxes condition, i.e
\begin{equation}
\label{ptr11}
\mu _K \left( {\nabla _{\mathcal D,\mu } u^{(i)}_h} \right)_{\left| {\left( {\bold x_{K^* } ,\bold x_K ,\bold x_\sigma  } \right)} \right.} .n^K_{[\bold x_\sigma, \bold x_{K^*}]}  + \mu _L \left( {\nabla _{\mathcal D,\mu } u^{(i)}_h} \right)_{\left| {\left( {\bold x_{K^* } ,\bold x_L ,\bold x_\sigma  } \right)} \right.} .n^L_{[\bold x_\sigma, \bold x_{K^*}]}  = 0,
\end{equation}
where $\mu_K$, $\mu_L$ are the average values of $\mu$ on $K$ and $L$.\\
Equation (\ref{ptr11}) leads to the following linear combination depended on $\{u^{(i)}_K, u^{(i)}_L,u^{(i)}_{K^*}  \}$
\begin{equation}
\label{ptr12}
u_{\sigma,K^*}^{(i)}  = \beta _K^{K^* ,\sigma } u^{(i)}_K  + \beta _L^{K^* ,\sigma } u^{(i)}_L  + \beta_{K^*}^{K^* ,\sigma } u^{(i)}_{K^*}, \quad \text{for each}~i = 1,2, 
\end{equation}
where the coefficients are written by
\begin{displaymath}
\beta _K^{K^* ,\sigma }  = {{\left( {\frac{{(n_{[\bold x_\sigma  ,\bold x_{K^* } ]}^K )^T \mu _K n_{[\bold x_\sigma  ,\bold x_{K^* } ]}^K }}{{2m_{(\bold x_{K^* } ,\bold x_{K,} \bold x_\sigma  )} }}} \right)} \mathord{\left/
 {\vphantom {{\left( {\frac{{(n_{[\bold x_\sigma  ,\bold x_{K^* } ]}^K )^T \mu _K n_{[\bold x_\sigma  ,\bold x_{K^* } ]}^K }}{{2m_{(\bold x_{K^* } ,\bold x_{K,} \bold x_\sigma  )} }}} \right)} {{\left( { - \frac{{(n_{[\bold x_\sigma  ,\bold x_{K^* } ]}^K )^T \mu _K n_{[\bold x_{K^* } ,\bold x_K ]} }}{{2m_{(\bold x_{K^* } ,\bold x_{K,} \bold x_\sigma  )} }} - \frac{{(n_{[\bold x_\sigma  ,\bold x_{K^* } ]}^L )^T \mu _L n_{[\bold x_{K^* } ,\bold x_L ]} }}{{2m_{(\bold x_{K^* } ,\bold x_{L,} \bold x_\sigma  )} }}} \right)}}}} \right.
 \kern-\nulldelimiterspace} {{\left( { - \frac{{(n_{[\bold x_\sigma  ,\bold x_{K^* } ]}^K )^T \mu _K n_{[\bold x_{K^* } ,\bold x_K ]} }}{{2m_{(\bold x_{K^* } ,\bold x_{K,} \bold x_\sigma  )} }} - \frac{{(n_{[\bold x_\sigma  ,\bold x_{K^* } ]}^L )^T \mu _L n_{[\bold x_{K^* } ,\bold x_L ]} }}{{2m_{(\bold x_{K^* } ,\bold x_{L,} \bold x_\sigma  )} }}} \right)}}},
\end{displaymath}
\begin{displaymath}
\beta _L^{K^* ,\sigma }  = {{\left( {\frac{{(n_{[\bold x_\sigma  ,\bold x_{K^* } ]}^L )^T \mu _L n_{[\bold x_\sigma  ,\bold x_{K^* } ]}^L }}{{2m_{(\bold x_{K^* } ,\bold x_{L,} \bold x_\sigma  )} }}} \right)} \mathord{\left/
 {\vphantom {{\left( {\frac{{(n_{[\bold x_\sigma  ,\bold x_{K^* } ]}^L )^T \mu _K n_{[\bold x_\sigma  ,\bold x_{K^* } ]}^L }}{{2m_{(\bold x_{K^* } ,\bold x_{L,} \bold x_\sigma  )} }}} \right)} {\left( { - \frac{{(n_{[\bold x_\sigma  ,\bold x_{K^* } ]}^K )^T \mu _K n_{[\bold x_{K^* } ,\bold x_K ]} }}{{2m_{(\bold x_{K^* } ,\bold x_{K,} \bold x_\sigma  )} }} - \frac{{(n_{[\bold x_\sigma  ,\bold x_{K^* } ]}^L )^T \mu _L n_{[\bold x_{K^* } ,\bold x_L ]} }}{{2m_{(\bold x_{K^* } ,\bold x_{L,} \bold x_\sigma  )} }}} \right)}}} \right.
 \kern-\nulldelimiterspace} {\left( { - \frac{{(n_{[\bold x_\sigma  ,\bold x_{K^* } ]}^K )^T \mu _K n_{[\bold x_{K^* } ,\bold x_K ]} }}{{2m_{(\bold x_{K^* } ,\bold x_{K,} \bold x_\sigma  )} }} - \frac{{(n_{[\bold x_\sigma  ,\bold x_{K^* } ]}^L )^T \mu _L n_{[\bold x_{K^* } ,\bold x_L ]} }}{{2m_{(\bold x_{K^* } ,\bold x_{L,} \bold x_\sigma  )} }}} \right)}},
\end{displaymath}
$\beta _{K^* }^{K^* ,\sigma }  = 1 - \beta _K^{K^* ,\sigma }  - \beta _L^{K^* ,\sigma } $.\\
\\
From Equation (\ref{ptr12}), the  unknown $u^{(i)}_{\sigma,K^*}$ is computed by $u^{(i)}_K,~u^{(i)}_{K^*}$ and $u^{(i)}_L$. Thus, the discrete gradient $\nabla _{\mathcal D,\mu } u_h^{(i)}$ on $(\bold x_K,\bold x_L,\bold x_{K^*})$ only depends on these three values.\\
\\
{\bf Hypothesis \ref{sec3}.1}: we assume  
\begin{equation}
\label{ptr13}
{\left( { - \frac{{(n_{[\bold x_\sigma  ,\bold x_{K^* } ]}^K )^T \mu _K n_{[\bold x_{K^* } ,\bold x_K ]} }}{{2m_{(\bold x_{K^* } ,\bold x_{K,} \bold x_\sigma  )} }} - \frac{{(n_{[\bold x_\sigma  ,\bold x_{K^* } ]}^L )^T \mu _L n_{[\bold x_{K^* } ,\bold x_L ]} }}{{2m_{(\bold x_{K^* } ,\bold x_{L,} \bold x_\sigma  )} }}} \right)} \ne 0.
\end{equation}
Note that if the mesh points $\bold x_K$ or $\bold x_L$ are moved slightly, the value of the right hand side in (\ref{ptr13}) is changed. This asserts that the hypothesis \ref{sec3}.1 is easy to be satisfied.\\
\newline
Using the definition of the discrete gradient $\nabla_\mathcal D u_h^{(i)}$, for each element of $\bold u_h$, we define the discrete divergence on the triangle $(\bold x_{K^*},\bold x_K,\bold x_L)$  by
\begin{equation}
\label{ptr14}
\text{div}_{\mathcal D,\mu}(\bold u_h) = \nabla_{\mathcal D,\mu} u_h^{(1)}.\bold e_1 + \nabla_{\mathcal D,\mu} u_h^{(2)}.\bold e_2
\end{equation}
with $\bold e_i$, $i=1,2$, the basis unit vector corresponding to the {\it i-th} coordinate.
\subsection{Discrete variational formulation}
\label{family6}
The existence and uniqueness of weak solution of the Stokes problem (\ref{ptr1}) was stated in the section \ref{sec1a}. However, in order to apply the pFECC scheme, we would like to implement another usual variational formulation for the problem (\ref{ptr1}), as follows: Find the velocity $\bold u \in (H^1_0(\Omega))^2$ and the pressure $p \in L^2_0(\Omega)$ such that
\begin{equation}
\label{weakptr2}
\left\{ \begin{array}{l}
\int\limits_\Omega  {\mu (\bold x)\nabla \bold u:\nabla \bold v\,d\bold x}  - \int\limits_\Omega  {\text{div} (\bold v)\,p\,dx}  = \int\limits_\Omega  {\bold f.\bold v\,d\bold x} ,\quad \bold v \in {(H_0^1(\Omega ))^2}\\
\int\limits_\Omega  {\text{div}(\bold u)\,q\,d\bold x}  = 0,\quad q \in L_0^2(\Omega ),
\end{array} \right.
\end{equation}
with the space
\begin{displaymath}
L^2_0(\Omega) = \{p \in L^2(\Omega) | \int\limits_\Omega{p~ \text{d} \bold x} = 0 \}.
\end{displaymath}
Applying the pFECC scheme into the velocity-pressure system (\ref{weakptr2}), we will look for the discrete velocity $\bold u_h \in \mathcal H_{\mathcal D}$, the discrete pressure $p_h \in \mathcal L_{\mathcal D}$ satisfying the following  problem
\begin{eqnarray}
\label{disweakptr2a}
&&\int\limits_\Omega  {\mu (\bold x)\nabla_{\mathcal D,\mu} \bold u_h:\nabla_{\mathcal D,\mu} \bold v_h\, \text{d}\bold x}  - \int\limits_\Omega  {\text{div}_{\mathcal D,\mu} (\bold v_h)\,p_h\,\text{d}\bold x}  = \int\limits_\Omega  {\bold f.P(\bold v_h)\,d\bold x} , \nonumber \\
&&\text{with} ~P(\bold v_h) = (P (v_h^{(1)}),P (v_h^{(2)})),~\text{for all}~\bold v_h=(v^{(1)},v^{(2)}) \in \mathcal H_{\mathcal D},
\end{eqnarray}
and
\begin{equation}
\label{disweakptr2b}
\int\limits_\Omega  {\text{div}_{\mathcal D,\mu}(\bold u_h)\,q_h\,\text{d}\bold x}  = -\lambda h\int\limits_\Omega {p_h q_h}~ \text{d} \bold x,\quad \text{for all}~q_h \in \mathcal L_{\mathcal D},
\end{equation}
where $\nabla_{\mathcal D,\mu} \bold u_h:\nabla_{\mathcal D,\mu} \bold v_h$ is defined by 
\begin{equation}
\label{ptr15}
\nabla_{\mathcal D,\mu} \bold u_h:\nabla_{\mathcal D,\mu} \bold v_h = \nabla_{\mathcal D,\mu} u_h^{(1)}.\nabla_{\mathcal D,\mu} v_h^{(1)} + \nabla_{\mathcal D,\mu} u_h^{(2)}.\nabla_{\mathcal D,\mu} v_h^{(2)}.
\end{equation}
\begin{equation}
\label{ptr15a}
\text{div}_{\mathcal D,\mu} \bold u_h = \partial _{\mathcal D,\mu }^{(1)}{u_h^{(1)}} + \partial _{\mathcal D,\mu }^{(2)}{u_h^{(2)}}
\end{equation}
where the discrete partial divergence ${\partial^{j} _{\mathcal D,\mu }}{u_h^{(i)}}$ corresponds to a discretization of the partial divergence $\frac{{\partial {u_h^{(i)}}}}{{\partial {x_j}}}$, defined by 
\begin{equation}
\label{ptr15b}
{\partial^{j} _{\mathcal D,\mu }}{u_h^{(i)}} = \nabla_{\mathcal D, \mu} u_h^{(i)}.\bold e_j \quad \text{for}~i,j \in \{1,2 \}.
\end{equation}
\subsection{The linear algebraic systems}
\label{family7}
Let us describe the three implementation steps to construct the system of linear equations depended on 
$\{\bold u_{K} \}_{K \in \mathcal M}$ and $\{p_{K^*}\}_{K^* \in \mathcal M^*}$, as follows: \\
\newline
{\em In the first step:} For each element $K^* \in \mathcal M^*$, a discrete test pressure function $q_h \in \mathcal L_{\mathcal D}$ is only equal to $1$ on $K^*$ and $0$ on $L^* \in \mathcal M^*/\{K^*\}$, Equation(\ref{disweakptr2b}) is stated by
 \begin{equation}
\label{ptr16}
\int\limits_{{K^*}} {\left( {{\nabla _{\mathcal D,\mu }}{u_h^{(1)}}.{\bold e_1} + {\nabla _{\mathcal D,\mu }}{u_h^{(2)}}.{\bold e_2}} \right) \text{d} \bold x}  =  - \lambda ~{h}~ {m}\left( {{K^*}} \right){p_{{K^*}}}
\end{equation}
Besides, Equation (\ref{disweakptr2a}) is computed with each value of a discrete test velocity function:
\begin{itemize}
\item $\bold v_h = (\{\bold v_{L}\}_{L \in \mathcal M},\{\bold v_{L^*}\}_{L^* \in \mathcal M^*})$  satisfies  $\bold {v}_{K^*} = (1, 0)$; $\bold v_{M} = (0, 0)$ for all $M \in \{\mathcal M \cup \mathcal M^*\} /\{K^*\}$, then this equation is rewritten as
\begin{equation}
\label{ptr17}
\int\limits_{{K^*}} {\left( {\mu (\bold x){\nabla _{\mathcal D,\mu }}{u_h^{(1)}}} \right).{\nabla _{\mathcal D,\mu }} v_h{^{(1)}}\text{d} \bold x}  - \int\limits_{{K^*}} {\left( {{\nabla _{\mathcal D,\mu }} v_h{\,^{(1)}}.{\bold e_1}} \right)p~\text{d} \bold x}  = \int\limits_{{K^*}} {{f_1}~\text{d} \bold x},
\end{equation}
\item $\bold w_h = (\{\bold w_{L}\}_{L \in \mathcal M},\{\bold w_{L^*}\}_{L^* \in \mathcal M^*})$  satisfies  $\bold {w}_{K^*} = (0, 1)$; $\bold w_{M} = (0, 0)$ for all $M \in \{\mathcal M \cup \mathcal M^*\} /\{K^*\}$, then it is equal to
\begin{equation}
\label{ptr18}
\int\limits_{{K^*}} {\left( {\mu (\bold x){\nabla _{\mathcal D,\mu }}{u_h^{(2)}}} \right).{\nabla _{\mathcal D,\mu }} v_h{^{(2)}}\text{d} \bold x}  - \int\limits_{{K^*}} {\left( {{\nabla _{\mathcal D,\mu }} v_h{^{(2)}}.{\bold e_2}} \right)p~\text{d} \bold x}  = \int\limits_{{K^*}} {{f_2}~\text{d} \bold x},
\end{equation}
\end{itemize}
We see that Equations (\ref{ptr16}), (\ref{ptr17}) and (\ref{ptr18}) can be represented as three linear combinations only depending on $\{\bold u_K\}_{K \in \mathcal M}$, $p_{K^*}$, $\bold u_{K^*} =(u^{(1)}_{K^*}), u^{(2)}_{K^*})$ and $\bold f$. These results help us compute the unknown $u^{(i)}_{K^*}$ by a linear combination $\Pi^{(i)}_{{K^*}}(\{u^{(i)}_K\}_{K \in \mathcal M}, \bold f) + \Psi^{(i)}(p_{K^*})$  with $i = 1, 2$. \\
\newline
{\bf Remark \ref{sec3}.1}: The coefficients of $u^{(i)}_K$ for all $K \in \mathcal M$, $i=1,2$, in the operator $\Pi^{(i)}_{{K^*}}$ are the same as those of 
the function $\Pi_{K^*}$, in the second step of \cite{PoHa}.\\
\newline 
{\em In the second step:} The unknowns $\bold u_{K^*}$, for all $K^* \in \mathcal M^*$, in the discrete gradient $\nabla_{\mathcal D,\mu} \bold u_h$ and the discrete divergence $\text{div}_{\mathcal D,\mu} \bold u_h$, are transformed into $\Pi^{(i)}_{\bold u_{K^*}}(\{u^{(i)}_K\}_{K \in \mathcal M}, \bold f) + \Psi^{(i)}(p_{K^*})$ and Ep.(\ref{ptr16}) with $i = 1, 2$.\\
\newline 
{\em In the last step:} For each element $K \in \mathcal M$, a test velocity function in Equation (\ref{disweakptr2a}) is taken into each following value:
\begin{itemize}
\item $\bold v = (\{\bold v_{K}\}_{K \in \mathcal M},\{\bold v_{K^*}\}_{K^* \in \mathcal M^*})$ has $\bold {v}_{K} = (1, 0)$; $\bold v_{M} = (0, 0)$ for all $M \in \{\mathcal M \cup \mathcal M^*\} /\{K\}$, the equation is computed by
\begin{equation}
\label{ptr19}
\int\limits_{{\Omega}} {\left( {\mu (\bold x){\nabla _{\mathcal D,\mu }}{u_h^{(1)}}} \right).{\nabla _{\mathcal D,\mu }} v_h{^{(1)}}\text{d} \bold x}  - \int\limits_{{\Omega}} {\left( {{\nabla _{\mathcal D,\mu }} v_h{\,^{(1)}}.{\bold e_1}} \right)p_h~\text{d} \bold x}  =  \int\limits_{\Omega} {{\bold f}(x).{\bold v_h}\,\text{d}{\bold x} },
\end{equation}
\item $\bold w = (\{\bold w_{K}\}_{K \in \mathcal M},\{\bold w_{K^*}\}_{K^* \in \mathcal M^*})$ has $\bold w_K = (0, 1)$ and $\bold w_{M} = (0, 0)$ for all $M \in \{\mathcal M \cup \mathcal M^*\} /\{K\}$, then the equation is equal to 
\begin{equation}
\label{ptr20}
\int\limits_{{\Omega}} {\left( {\mu (\bold x){\nabla _{\mathcal D,\mu }}{u_h^{(2)}}} \right).{\nabla _{\mathcal D,\mu }} w_h{^{(2)}}\text{d} \bold x}  - \int\limits_{{\Omega}} {\left( {{\nabla _{\mathcal D,\mu }} w_h{\,^{(2)}}.{\bold e_2}} \right)p_h~\text{d} \bold x}  = \int\limits_{\Omega} {{\bold f}(x).{\bold w_h}\,\text{d}{\bold x} },
\end{equation}
\end{itemize}
{\bf Remark \ref{sec3}.2:} After the second step,  the discrete gradient $\nabla_{\mathcal D, \mu}\bold u_h$ and divergence $\text{div}_{\mathcal D, \mu} \bold u_h$ are independent on the unknowns $\{\bold u_{K^*}\}_{K^* \in \mathcal M^*}$. Therefore, in two equations (\ref{ptr19}) and (\ref{ptr20}), there are not the unknowns $\{\bold u_{K^*}\}_{K^* \in \mathcal M^*}$. Additionally, for each $i = 1,2$, we have  $supp\{\nabla_{\mathcal D,\mu}(u_h^{i})\}$ belonging to  $\bigcup\limits_{{K^*} \in \mathcal M_K^*} {{K^*}}$ with 
\begin{displaymath}
\mathcal M_K^* = \{ K^* \in \mathcal M^* ~|~ K \cap K^* \ne \emptyset \},
\end{displaymath}
which indicates that the stiffness matrix $\bold A$ in (\ref{ptr21})    is sparse\\
\newline
From the above three steps, the stiffness matrix associated to the pFECC scheme for the Stokes problem is generated by
\begin{equation}
\label{ptr21}
\underbrace {\left( {\begin{array}{*{20}{c}}
{{{\bold{B}}}}& \bold C\\
{{{\bold{C}}}}&{{(-\lambda~h~ {m}(K^*)\bf{Id})}}
\end{array}} \right)}_{\bold A}\left( \begin{array}{l}{\{ {\bold u_K}\} _{K \in \mathcal M}}\\
{\{{p_{{K^*}}}\} _{{K^*} \in {\mathcal M^*}}}
\end{array} \right) = {\bold{F}}
\end{equation}
Thanks to Remark \ref{sec3}.1, the matrix $\bold B$ are positive definite and symmetric, which is proven in Lemma $3.2$ of \cite{PoHa}. Hence, the matrix $\bold A$ has the inverse matrix, there then exists the unique solution of the system (\ref{ptr21}).
\section{Consistency and stability of the pFECC method}
\label{sec4}
In this section, we will study the consistency properties of the discrete gradient and the discrete divergence. These results will be necessary to prove the convergence. Let us firstly consider to the divergence operator.\\
\newline
{\bf Lemma} \ref{sec4}.1 {\bf (Consistency of the discrete divergence)}\\
\textit{Under geometrical conditions for meshes \ref{family4} and assumptions of Lemma $5.1$ in \cite{PoHa}, then, there exists the positive constant $C_5$, such that, for all $ \overline{ u} \in H^2(\Omega)\cap H^1_0(\Omega)$, and for each $K^* \in \mathcal M^*$,}
\begin{equation}
\label{ptr22}
\left| {\int\limits_{K^*} {{\partial^{(i)}_{\mathcal D,\mu }}\overline{u}_h(\bold x)\,\text{d}\bold x}  - \int\limits_{{\partial K^*}} {\overline{u}(x){\bold e_i}.n(\bold x)\,\text{d}\bold x}} \right| \le {C_4}h^2 \left[ |{\overline{u}}{|_{{H^2}({K^*})}}+  ||\nabla \overline u||_{L^2(K^*)} \right]
\end{equation}
with $\overline u_h = (\{\overline u(\bold x_{K})\}_{K \in \mathcal M},\{\overline u(\bold x_{K^*})\}_{K^* \in \mathcal M^*})$ and $\partial^{(i)}_{\mathcal D,\mu }\overline{u}_h = \nabla_{\mathcal D,\mu}\overline{u}_h . \bold e_i$.\\
\newline
\textit{Proof.}\\
For any $K^* \in \mathcal M^*$, we have
\begin{equation}
\label{ptr23}
{\int\limits_{K^*} {\partial^{(i)}_{\mathcal D,\mu }}{\overline {u}_h \text{d}\bold x} }
 = \sum\limits_{{L^{**}} \in \mathcal M_{{K^*}}^{**}} ~ {\int\limits_{{L^{**}}} {\partial^{(i)}_{\mathcal D,\mu }\overline {u}_h\,\text{d}\bold x} }
\end{equation}
with $\mathcal M_{{K^*}}^{**} = \{ L^{**} \in \mathcal M^{**} ~|~ L^{**} \subset K^*\}$.\\
\newline
Let us consider on any element $K^{**} \in  \mathcal M^{**}_{K^*}$, it is seen as a triangle $(\bold x_{K^*},\bold x_K, \bold x_L)$ having three vertices $\bold x_K$, $\bold x_L$ with $K, L \in \mathcal M$ and $\bold x_{K^*}$ with $K^* \in \mathcal M^*$.\\
\newline
On this triangle, we compute 
\begin{equation}
\label{ptr23a}
\int\limits_{\left( {{\bold x_K},{\bold x_L},{\bold x_{{K^*}}}} \right)} {\partial _{\mathcal D,\mu }^{(i)}\overline u_h(\bold x)d\bold x}  = \int\limits_{\left( {{\bold x_K},{\bold x_{{K^*}}},{\bold x_\sigma }} \right)} {{\nabla _{\mathcal D,\mu }}\overline u_h.{\bold e_i}\text{d}\bold x}  + \int\limits_{\left( {{\bold x_L},{\bold x_{{K^*}}},{\bold x_\sigma }} \right)} {{\nabla _{\mathcal D,\mu }}\overline u_h .{\bold e_i}\text{d}\bold x} 
\end{equation}
with
\begin{equation}
\label{ptr23b}
\int\limits_{\left( {{\bold x_K},{\bold x_{{K^*}}},{\bold x_\sigma }} \right)} {{\nabla _{\mathcal D,\mu }}\bar u.{\bold e_i}\text{d}\bold x} = -\frac{{1}}{{2}}
\left\{ \begin{array}{l}
\vspace{0.5cm}
\left( {{n_{[{\bold x_K},{\bold x_\sigma }]}} + \beta _{{K^*}}^{{K^*},\sigma }{n_{\left[ {{\bold x_{{K}}},{\bold x_{K^*} }} \right]}}} \right).{\bold e_i}\,{\overline u(\bold x_{{K^*}})} \\
\vspace{0.5cm}
+ \beta _L^{{K^*},\sigma }{n_{\left[ {{\bold x_K},{\bold x_{{K^*}}}} \right]}}.{\bold e_i}\,{\bar u(\bold x_L)}\\
+ \left( {{n^{(K,K^*,\sigma)}_{[{\bold x_{K^*}},{\bold x_\sigma }]}} + \beta _K^{{K^*},\sigma }{n_{\left[ {{\bold x_{{K^*}}},{\bold x_K }} \right]}}} \right).{\bold e_i}\,{\bar u(\bold x_K)}
\end{array} \right\}
\end{equation}
and 
\begin{equation}
\label{ptr23c}
\int\limits_{\left( {{\bold x_L},{\bold x_{{K^*}}},{\bold x_\sigma }} \right)} {{\nabla _{\mathcal D,\mu }}\bar u.{e_i}d\bold x} = -\frac{{1}}{{2}}
\left\{ \begin{array}{l}
\vspace{0.5cm}
\left( {{n_{[{\bold x_L},{\bold x_\sigma }]}} + \beta _{{K^*}}^{{K^*},\sigma }{n_{\left[ {{\bold x_{{L}}},{\bold x_{K^*} }} \right]}}} \right).{e_i}\,{\bar u(\bold x_{{K^*}})}\\
\vspace{0.5cm}
+ \beta _K^{{K^*},\sigma }{n_{\left[ {{\bold x_L},{\bold x_{{K^*}}}} \right]}}.{e_i}\,{\bar u(\bold x_K)}\\
+ \left( {{n^{(L,K^*,\sigma)}_{[{\bold x_{K^*}},{\bold x_\sigma }]}} + \beta _L^{{K^*},\sigma }{n_{\left[ {{\bold x_{L}},{\bold x_{K^*} }} \right]}}} \right).{e_i}\,{\bar u(\bold x_L)}
\end{array} \right\}
\end{equation}
In Equations (\ref{ptr23a})-(\ref{ptr23c}), the vectors $n_{[\cdot,\cdot]}$ which are outward normal vectors, their lengths are equal to the corresponding edge lengths. The three coefficients $\beta _{{K}}^{{K^*},\sigma }$, $\beta _{{L}}^{{K^*},\sigma }$ and $\beta _{{K^*}}^{{K^*},\sigma }$ are defined in Equation (\ref{ptr12}). We also define $\sigma = [\bold x_K,\bold x_L]$ and the common edge $K|L$.\\
\newline
Similarly, we transform the integral of $\overline{ u}$ on $\left( {{\bold x_K},{\bold x_L},{\bold x_{{K^*}}}} \right)$ into two integrals on two sub-triangles $\left( {{\bold x_K},{\bold x_{{K^*}}},{\bold x_\sigma }} \right)$ and $\left( {{\bold x_L},{\bold x_{{K^*}}},{\bold x_\sigma }} \right)$
\begin{equation}
\label{ptr24}
\int\limits_{\left( {{\bold x_K},{\bold x_L},{\bold x_{{K^*}}}} \right)} {\frac{{\partial \bar u}}{{\partial { x_i}}}} (\bold x)\text{d} \bold x = \int\limits_{\left( {{\bold x_K},{\bold x_{{K^*}}},{\bold x_\sigma }} \right)} {\frac{{\partial \bar u}}{{\partial { x_i}}}} (\bold x)\text{d}\bold x + \int\limits_{\left( {{\bold x_L},{\bold x_{{K^*}}},{\bold x_\sigma }} \right)} {\frac{{\partial \bar u}}{{\partial {x_i}}}} (\bold x)\text{d}\bold x,
\end{equation}
In the above equation, two integrals of the right hand side are computed by
\begin{eqnarray}
\label{ptr25}
\int\limits_{\left( {{\bold x_K},{\bold x_{{K^*}}},{\bold x_\sigma }} \right)} {\frac{{\partial \bar u}}{{\partial { x_i}}}} (\bold x)\text{d}\bold x &=& \int\limits_{\left[ {{\bold x_K},{\bold x_{{K^*}}}} \right]} {\bar u\,({\bold e_i}.\bar n_{\left[ {{\bold x_K},{\bold x_{{K^*}}}} \right]})} \text{d}\gamma (\bold x) + \int\limits_{[{\bold x_K},{\bold x_\sigma }]} {\bar u\,({\bold e_i}.\bar n_{[{\bold x_K},{\bold x_\sigma }]})} \text{d}\gamma (\bold x), \nonumber \\
&+& \int\limits_{[{\bold x_{K^*}},{\bold x_\sigma }]} {\bar u\,({\bold e_i}.\bar n^{(K,K^*,\sigma)}_{[{\bold x_{K^*}},{\bold x_\sigma }]})} \text{d}\gamma (\bold x).
\end{eqnarray}
and
\begin{eqnarray}
\label{ptr26}
\int\limits_{\left( {{\bold x_L},{\bold x_{{K^*}}},{\bold x_\sigma }} \right)} {\frac{{\partial \bar u}}{{\partial {x_i}}}} (\bold x)\text{d}\bold x &=&  \int\limits_{\left[ {{\bold x_L},{\bold x_{{K^*}}}} \right]} {\bar u\,({\bold e_i}.\bar n_{[{\bold x_L},{\bold x_{K^*} }]})} \text{d}\gamma (\bold x) + \int\limits_{[{\bold x_L},{\bold x_\sigma }]} {\bar u\,({\bold e_i}.\bar n_{[{\bold x_L},{\bold x_\sigma }]})} \text{d}\gamma (\bold x) \nonumber \\
&+& \int\limits_{[{\bold x_{K^*}},{\bold x_\sigma }]} {\bar u\,({\bold e_i}.\bar n^{(K,K^*,\sigma)}_{[{\bold {x}_{K^*}},{\bold x_\sigma }]})} \text{d}\gamma (\bold x),
\end{eqnarray}
where the vectors $\bar n_{[.,.]}$ are outward normal unit vectors 
of the considered triangle, and $\bar n^{(K,K^*,\sigma)}_{[{\bold x_{K^*}},{\bold x_\sigma }]} + \bar n^{(L,K^*,\sigma)}_{[{\bold x_{K^*}},{\bold x_\sigma }]} = \bold 0$. Besides, we have the following relationship between $\overline n_{[\cdot,\cdot]}$ and $n_{[\cdot,\cdot]}$
\begin{eqnarray}
\label{ptr27}
&&n_{[{\bold x_L},{\bold x_\sigma }]} = {m}([{\bold x_L},{\bold x_\sigma }]) \bar n_{[{\bold x_L},{\bold x_\sigma }]},~n_{[{\bold x_K},{\bold x_\sigma }]} = {m}([{\bold x_K},{\bold x_\sigma }]) \bar n_{[{\bold x_K},{\bold x_\sigma }]},\nonumber \\
&&n_{[{\bold x_L},{\bold x_{K^*} }]} = {m}([{\bold x_L},{\bold x_{K^*} }]) \bar n_{[{\bold x_L},{\bold x_{K^*} }]},~n_{[{\bold x_K},{\bold x_{K^*} }]} = {m}([{\bold x_K},{\bold x_{K^*} }]) \bar n_{[{\bold x_K},{\bold x_{K^*} }]}
\end{eqnarray}
From Equations (\ref{ptr23a})-(\ref{ptr27}), they lead
\begin{eqnarray}
\label{ptr28}
&&\int\limits_{\left( {{\bold x_K},{\bold x_L},{\bold x_{{K^*}}}} \right)} {\partial _{D,\mu }^{(i)}\overline u_h \text{d}\bold x}  - \int\limits_{\left( {{\bold x_K},{\bold x_L},{\bold x_{{K^*}}}} \right)} {\frac{{\partial \bar u}}{{\partial {x_i}}}} (\bold x)\text{d}\bold x = \nonumber \\
\vspace{0.7cm}
&=& \frac{1}{2}\left( {{n_{\left[ {{\bold x_{{K^*}}},{\bold x_K}} \right]}}.{{\bf{e}}_i}} \right)\left( {\bar u({\bold x_{{K^*}}}) + \bar u({\bold x_K}) - \frac{2}{{{{m}}\left[ {{\bold x_{{K^*}}},{\bold x_K}} \right]}}\int\limits_{\left[ {{\bold x_{{K^*}}},{\bold x_K}} \right]} {\overline u(\bold x)\text{d}\gamma (\bold x)} } \right) \nonumber \\
\vspace{0.7cm}
&+& \frac{1}{2}\left( {{n_{\left[ {{\bold x_{{K^*}}},{\bold x_L}} \right]}}.{{\bf{e}}_i}} \right)\left( {\bar u({\bold x_{{K^*}}}) + \bar u({\bold x_L}) - \frac{2}{{{{m}}\left[ {{\bold x_{{K^*}}},{\bold x_L}} \right]}}\int\limits_{\left[ {{\bold x_{{K^*}}},{\bold x_L}} \right]} {\bar u(\bold x)\text{d}\gamma (\bold x)} } \right)\nonumber \\
\vspace{0.7cm}
&+& \frac{1}{2}\left( {{n_{\left[ {{\bold x_K},{\bold x_L}} \right]}}.{{\bf{e}}_i}} \right)\left( {\bar u({\bold x_K}) + \overline u({\bold x_L}) - \frac{2}{{{{m}}\left[ {{\bold x_K},{\bold x_L}} \right]}}\int\limits_{\left[ {{\bold x_K},{\bold x_L}} \right]} {\overline u(\bold x)\text{d}\gamma (\bold x)} } \right)\nonumber \\
\vspace{0.7cm}
&+& \frac{1}{2}\left( {{n_{\left[ {{\bold x_K},{\bold x_L}} \right]}}.{{\bf{e}}_i}} \right)    \underbrace{\left[ {\bar u({\bold x_K})
\left(\frac{{{m}[{\bold x_L},{\bold x_\sigma }]}}{{{{m}}\left[ {{\bold x_K},{\bold x_L}} \right]}} - \beta _{{K}}^{{K^*},\sigma }\right) 
+ \bar u({\bold x_L})
\left( \frac{{m[{\bold x_K},{\bold x_\sigma }]}}{{{{m}}\left[ {{\bold x_K},{\bold x_L}} \right]}} - \beta _{{L}}^{{K^*},\sigma }\right) - \beta _{{K^*}}^{{K^*},\sigma }} \bar u(\bold x_{K^*}) \right]}
_{{[\bar u({\bold x_K})-\bar u(\bold x_{K^*})]
\left(\frac{{{m}[{\bold x_L},{\bold x_\sigma }]}}{{{{m}}\left[ {{\bold x_K},{\bold x_L}} \right]}} - \beta _{{K}}^{{K^*},\sigma }\right) 
+ [\bar u({\bold x_L})-\bar u(\bold x_{K^*})]
\left( \frac{{m[{\bold x_K},{\bold x_\sigma }]}}{{{{m}}\left[ {{\bold x_K},{\bold x_L}} \right]}} - \beta _{{L}}^{{K^*},\sigma }\right)} }. \nonumber \\
\end{eqnarray}
because of $\beta _{{K}}^{{K^*},\sigma } + \beta _{{L}}^{{K^*},\sigma }+\beta _{{K^*}}^{{K^*},\sigma } = 1$ and $\frac{{{m}[{\bold x_K},{\bold x_\sigma }]}}{{{{m}}\left[ {{\bold x_K},{\bold x_L}} \right]}} + \frac{{{m}[{\bold x_L},{\bold x_\sigma }]}}{{{{m}}\left[ {{\bold x_K},{\bold x_L}} \right]}} = 1$.\\
\newline
Note that we have the useful results for Equation (\ref{ptr28})
\begin{equation}
\label{ptr29}
\mathop {\lim }\limits_{h \to 0} 
\left(\frac{{{m}[{\bold x_L},{\bold x_\sigma }]}}{{{{m}}\left[ {{\bold x_K},{\bold x_L}} \right]}} - \beta _{{K}}^{{K^*},\sigma }\right)  = \mathop {\lim }\limits_{h \to 0} \left( \frac{{m[{\bold x_K},{\bold x_\sigma }]}}{{{{m}}\left[ {{\bold x_K},{\bold x_L}} \right]}} - \beta _{{L}}^{{K^*},\sigma }\right) = \mathop {\lim }\limits_{h \to 0} \beta _{{K^*}}^{{K^*},\sigma } = 0
\end{equation}
which is proven by Lemma $5.1$ in \cite{PoHa}, and
\begin{eqnarray}
\label{ptr29a}
|\overline u(\bold x_K) - \overline u(\bold x_{K^*})| \le C_{\ref{ptr29a}} ||\nabla \overline u||_{L^2(K^*)}\quad \text{for all}~K \in \mathcal M, K^* \in \mathcal M^*,
\end{eqnarray}
it is shown by Theorem $9.12$ (Morrey) in \cite{Bre}.\\
\newline
For the other computations of Eq.(\ref{ptr28}), let us give another triangular element $(\bold x_{K^*},\bold x_{K},\bold x_{M}) \in \mathcal M^{**}_{K^*}$. Two triangles $(\bold x_{K^*},\bold x_{K},\bold x_{L})$ and $(\bold x_{K^*},\bold x_{K},\bold x_{M})$ have a common edge $\left[ {{x_{K^*}},{x_{K}}} \right]$, so we should rewrite the vector $n_{\left[ {{\bold x_{K^*}},{\bold x_{K}}} \right]}$
by $n^{(K,M,K^*)}_{\left[ {\bold x_{K^*},\bold x_{K}}\right]}$. This help us distinguish the normal outward vector $n^{(K,L,K^*)}_{\left[ {\bold x_{K^*},\bold x_K} \right]}$ of $(\bold x_{K},\bold x_{L},\bold x_{K^*})$ and $n^{(K,M,K^*)}_{\left[ {\bold x_{K^*},\bold x_{K}} \right]}$ of $(\bold x_{K},\bold x_{M},\bold x_{K^*})$. Beside, we have the important property of the two normal vectors 
\begin{displaymath}
n^{(K,L,K^*)}_{\left[ {\bold x_{K^*},\bold x_{K}} \right]} + n^{(K,M,K^*)}_{\left[ {\bold x_{K^*},\bold x_{K}} \right]} = 0
\end{displaymath}
which also appears in the other edge having a common vertex $\bold x_{K^*}$.\\
\newline
Combining the above property, (\ref{ptr23}), (\ref{ptr28}) and (\ref{ptr29}), we can estimate the right hand size of Equation(\ref{ptr22}), as follows:
\begin{eqnarray}
\label{ptr30}
&& \left|\int\limits_{{K^*}} {\partial^{(i)}_{{\mathcal D,\mu }}\overline u_h \,\text{d}\bold x}  - \int\limits_{\partial {K^*}} {\bar u(\bold x){\bold e_i}.n(\bold x)\,d\bold x} \right| = \nonumber \\
&\le&  {\frac{1}{2}\sum\limits_{\sigma  \in \mathcal E_{{K^*}}^*} {\left|\left( {n_\sigma ^{{K^*}}.{e_i}} \right)\right| \left|\left( {\bar u({\bold x^\sigma_1}) + \bar u({\bold x^\sigma_2}) - \frac{2}{{m(\sigma )}}\int\limits_\sigma  {\bar u(\bold x)\,d\bold x} } \right)+ C_5 h   ||\nabla \overline u||_{L^2(K^*)}\right|}} ,\nonumber \\
\end{eqnarray}
where the points $\bold x^\sigma_1,~\bold x^\sigma_2$ are  two vertices of an edge $\sigma$, $\text{m}(\sigma)$ is its length, and $n^{K^*}_{\sigma}$ denotes the outward normal vector of $K^*$ at $\sigma$. \\
In Equation (\ref{ptr30}), we remain estimate
\begin{equation}
\label{ptr31}
{\psi _\sigma }(\bar u) = \bar u({\bold x^\sigma_1}) + \bar u({\bold x^\sigma_2}) - \frac{2}{{m(\sigma )}}\int\limits_\sigma  {\bar u(\bold x)\,d\gamma (\bold x)} 
\end{equation}
Let $T$ be the reference triangle with the three vertices $\bold x^T_1(0,0)$,  $\bold x^T_2(1,0)$, $\bold x^T_3(0,1)$, and we put $\theta$ be the affine mapping from $T_{K^*,\sigma} = (\bold x^\sigma_1,\bold x^\sigma_2,\bold x_{K^*})$ to $T$ such that $\theta(\bold x^\sigma_1) = \bold x^T_1$, $\theta(\bold x^\sigma_2) = \bold x^T_2$ and $\theta(\bold x_{K^*}) = \bold x^T_3$.\\
On the reference triangle $T$, the operator $\hat \psi$, which is defined in $(H^2(T))'$, satisfies the operator  $\hat \psi(\hat u) = \psi_\sigma(\hat u \circ \theta)$ and
\begin{equation}
\label{ptr32}
{\hat \psi _\sigma }(\hat u) = \hat u({\bold x^T_1}) + \hat u({\bold x^T_2}) - \frac{2}{{m\left( {\left[ {{\bold x^T_1},{\bold x^T_2}} \right]} \right)}}\int\limits_{\left[ {{\bold x^T_1},{\bold x^T_2}} \right]} {\hat u(\bold x)\,d\gamma(\bold x)}
\end{equation}
for all $\hat u \in H^2(T)$.\\
Thanks to the Bramble-Hilbert Lemma, we get the following estimation 
\begin{equation}
\label{ptr33}
|\hat \psi _{\sigma }(\hat u)| \le C_6 |\hat u|_{H^2(T)},
\end{equation}
where there exist the positive constant be independence with anything.\\
Using Inequality (\ref{ptr33}) to estimate ${\psi _\sigma }(\bar u)$, to this aim, we choose $\hat u = \bar u \circ \theta^{-1} \in H^2(T)$. It implies $\Psi_{\sigma}(\bar u) = \hat \Psi_{\sigma}(\hat u)$ and
\begin{equation}
\label{ptr34}
|\Psi_{\sigma}(\bar u)|  \le C_7 |\hat u|_{H^2(T)}.
\end{equation}
In order to complete the estimation (\ref{ptr34}), we use the theorems $3.1.2$, $3.1.3$ in \cite{Ciar} to give the following results
\begin{equation}
\label{ptr35}
|\hat u|_{H^2(T)} \le ||\theta^{-1}||^2 \left(\frac{{m}(T)}{{m}(T_{K^*,\sigma})}\right)^{1/2} |\bar u|_{H^2(T_{K^*,\sigma})},
\end{equation}
and $||\theta^{-1}|| \le \frac{\text{diam}(T_{K^*,\sigma}))}{\hat \rho}$, where $\hat \rho$ is a diameter of a inscribed circle in $T$.\\
Applying Inequality (\ref{ptr35}), $m(T) = 1$, $\hat \rho = 1$ and  the geometrical condition (\ref{ptr10a}) to Inequality (\ref{ptr34}), then it yields
\begin{equation}
\label{ptr36}
|\Psi_{\sigma}(\bar u)| \le C_7 \text{diam}(T_{K^*,\sigma }) |\bar u|_{H^2(T_{K^*,\sigma})} \le C_7 h  |\bar u|_{H^2(T_{K^*,\sigma})}.
\end{equation}
From the estimation (\ref{ptr36}) for each edge $\sigma \in \mathcal E^*_{K^*}$ and the geometrical condition (\ref{ptr7}), the right hand side of Inequality (\ref{ptr30}) is less than $\left( C_7 {h} ||u||_{H^2(K^*)} + C_5 h ||\nabla \overline u||_{L^2(K^*)} \right)$.\\
Moreover, let us any $K^* \in \mathcal M^*$, and $\sigma \in \mathcal E^*_{K^*}$, $\sigma$ is also an edge of the triangular mesh $\mathcal M^{**}$, it leads 
\begin{equation}
\label{ptr36a} 
|(n^{K^*}_{\sigma}.\bold e_i)| \le h \quad\text{for all}~i=1,2
\end{equation}
By the inequalities (\ref{ptr30}), (\ref{ptr36}) and (\ref{ptr36a}), we have 
\begin{eqnarray}
\label{ptr37a}
\left|\int\limits_{{K^*}} {\text{div}_{\mathcal D,\mu}\overline u_h\,\text{d}\bold x}  - \int\limits_{\partial {K^*}} {\bar u(\bold x){\bold e_i}.n(\bold x)\,\text{d}\bold x} \right| & \le & \frac{C_7 h^2}{2}  \sum\limits_{\sigma  \in \mathcal E_{{K^*}}^*} {|\bar u|_{H^2(T_{K^*,\sigma})}} \nonumber \\
 &+& \frac{C_5 h^2}{2}   ||\nabla \overline u||_{L^2(K^*)} ,\nonumber \\
\end{eqnarray}
Remark that Inequality (\ref{ptr37a}) only require $\bar u \in H^2(T) \cap H^1_0(\Omega)$ for all triangle $T \in \mathcal M^{**}$ with $\bigcup\limits_{T \in {\mathcal M^{**}}} \overline T  = \overline \Omega$.\\
Together $\bar u \in H^2(\Omega) \cap H_0^1(\Omega)$ , it follows
\begin{eqnarray}
\label{ptr37}
\left|\int\limits_{{K^*}} {\text{div}_{\mathcal D,\mu}\overline u_h\,d\bold x}  - \int\limits_{\partial {K^*}} {\bar u(\bold x){e_i}.n(\bold x)\,d\bold x} \right| & \le & h^2 C_4 \left[ |u|_{H^2(K^*)} + ||\nabla \overline u||_{L^2(K^*)} \right], \nonumber \\
\end{eqnarray}
where a positive constant $C_4 = \frac{1}{2}\max\{C_5,C_7\}$ is independent on $h$, and $h$ is small enough. $\square$\\
\newline
{\bf Lemma} \ref{sec4}.2 {\bf (Stability of the scheme)}\\
\textit{Under the geometrical conditions for meshes are satisfied. Then, there exist a positive constant (independent on $h$), such that}
\begin{equation}
\label{stab1}
\mathop {\sup }\limits_{\bold{v_h} \in {\mathcal H_{{\mathcal D}}}\hfill\atop
\bold{v_h} \ne 0\hfill} \frac{{\int\limits_\Omega  {\left( {\text{div}_{{\mathcal D,\mu }}{\bold v_h}} \right){q_h}\,{\text{d}}\bold x } }}{{||P_1(\bold  v_h)||_{(H^1(\Omega))^2}}} \ge C_{\ref{stab1}}||{q_h}|{|_{{L^2}(\Omega )}}
\end{equation}
\textit{for all $q_h \in \mathcal L_{\mathcal D}$, in which $P_1(\bold v_h)=\left(P_1(v^{(1)}_h),P_1(v^{(2)}_h) \right)$ is defined the traditional interpolation, constructed on $\mathcal M^{**}$, the basis Lagrange polynomials having the degree $1$ and $\bold v_h \in \mathcal H_{\mathcal D}$}.\\
\newline
\textit{Proof}\\
For the left hand size of (\ref{stab1}), we have
\begin{eqnarray}
\label{stab2}
&&\mathop {\sup }\limits_{\bold{v}_h \in {\mathcal H_{{\mathcal D}}}\hfill\atop
{\bold v_h} \ne 0\hfill}\int\limits_\Omega  { \frac{\left({\text{div}_{{\mathcal D,\mu }}{\bold v_h}} \right)}{||P_1(\bold  v_h)||_{(H^1(\Omega))^2}} {q_h}\,{\text{d}}\bold x } \nonumber \\
& =& \mathop {\sup }\limits_{{\bold v_h} \in {\mathcal H_{{\mathcal D}}}\hfill\atop
{\bold v_h} \ne 0\hfill} \left \{ \int\limits_\Omega  {\frac{\left[ {\text{div}_{{\mathcal D,\mu }}{\bold v_h}} - \text{div}(P_1(\bold v_h)\right]}{||P_1(\bold  v_h)||_{(H^1(\Omega))^2}}{q_h}\,{\text{d}}\bold x } +  \int\limits_\Omega  { \frac{\text{div}(P_1(\bold v_h))}{||P_1(\bold  v_h)||_{(H^1(\Omega))^2}} {q_h}\,{\text{d}}\bold x } \right \} \nonumber \\
&\nonumber\\
&\ge& \mathop {\sup }\limits_{{\bold v_h} \in {\mathcal H_{{\mathcal D}}}\hfill\atop
{\bold v_h} \ne 0\hfill} \left [ \int\limits_\Omega  { \frac{\text{div}(P_1(\bold v_h))}{||P_1(\bold  v_h)||_{(H^1(\Omega))^2}} {q_h}\,{\text{d}}\bold x } \right ] -  \mathop {\sup }\limits_{{\bold v_h} \in {\mathcal H_{{\mathcal D}}}\hfill\atop
\bold{v_h} \ne 0\hfill} \left| \int\limits_\Omega  {\frac{\left[ {\text{div}(P_1(\bold v_h) - \text{div}_{{\mathcal D,\mu }}{\bold v_h}} \right]}{||P_1(\bold  v_h)||_{(H^1(\Omega))^2}}{q_h}\,{\text{d}}\bold x } \right |  \nonumber \\
\end{eqnarray}
where let any $\bold v_h \in \mathcal H_{\mathcal D}$, the function $P_1(\bold v_h)$ is a linear combination constructed by Lagrange polynomials of degree one on $\mathcal M^{**}$ and values of all element of $\bold v_h$.\\
\newline
To estimate $\int\limits_\Omega  { \frac{\text{div}(P_1(\bold v_h))}{||P_1(\bold  v_h)||_{(H^1(\Omega))^2}} {q_h}\,{\text{d}}\bold x }$, we recall, in the construction of the primal $\mathcal M$ and dual $\mathcal M^*$ meshes, we see that each element of $\mathcal M^*$ containing at most fixed number, $C_1~ - $ this is the condition (\ref{ptr7}), triangles of $\mathcal M$. From this property, each element of $\mathcal M^*$ forms a disjoint polygonal "macroelement". Moreover, according to Definition \ref{sec3}.2, each $p_h \in \mathcal L_{\mathcal D_h}$ is piecewise constant on each "macroelement" $K^* \in \mathcal M^*$. We then apply to the macroelement technique in \cite{CrRa}, \cite{Sten1} and Theorem $3$ in \cite{BrePit}. This leads the stability property is satisfied by using the Fortin's trick \cite{BoBrFo} for checking the inf-sup condition, i.e, 
\begin{equation}
\label{stab3}
\mathop {\sup }\limits_{\bold{v}_h \in {\mathcal H_{{\mathcal D}}}\hfill\atop
~\bold{v_h} \ne 0\hfill} \left[ \int\limits_\Omega  { \frac{\text{div}(P_1(\bold v_h))}{||P_1(\bold  v_h)||_{(H^1(\Omega))^2}} {q_h}\,{\text{d}}\bold x } \right ] \ge 
\mathop {\sup }\limits_{\bold{w} \in {(\mathbb V_h(\mathcal M^{**}))^2}\hfill\atop
\quad~ \bold{w} \ne 0\hfill} \left[ \int\limits_\Omega  { \frac{\text{div}(\bold w)}{||\bold w||_{{(H^1(\Omega))^2}}} {q_h}\,{\text{d}}\bold x } \right ] 
\ge \eta ||q_h||_{L^2(\Omega)} ,
\end{equation}
where the positive constant $\eta$ is independent on $h$, $\mathbb V_h(\mathcal M^{**})$ is the finite element space of the standard finite element method on the triangulation $\mathcal M^{**}$. Remark that 
\begin{displaymath}
(\mathbb V_h(\mathcal M^{**}))^2 \subseteq \left\{ {{P_1}({\bold v_h})|\forall {\bold v_h} \in {\mathcal H_{\mathcal D}}} \right\}
\end{displaymath}
Next, we estimate the following integral
\begin{eqnarray}
\label{stab4}
&&\int\limits_\Omega  {q_h \frac{{\left[ {\text{div}({P_1}({v_h})) - \text{div}_{{D,\mu }}{v_h} } \right]}}{{||P_1(\bold  v_h)||_{(H^1(\Omega))^2}}}} {\rm{d}}{\bf{x}} =  \sum\limits_{{K^*} \in {\mathcal M^*}} {\int\limits_{{K^*}} {{{\rm{q}}_h}\sum\limits_{i = 1}^2 {\frac{{\left[ {\partial _{D,\mu }^{(i)}v_h^{(i)} - {\partial ^{(i)}}{P_1}\left( {v_h^{(i)}} \right)} \right]}}{{||P_1(\bold  v_h)||_{(H^1(\Omega))^2}}}} {\rm{d}}{\bf{x}}} } \nonumber \\
&\nonumber \\
&\mathop {\le}\limits_{(\ref{ptr37a})}  &  \sum\limits_{{K^*} \in {\mathcal M^*}}  |q_{K^*}|{ { \sum\limits_{i = 1}^2 {h^2_{\mathcal D^{**}}\frac{\frac{C_{7}}{2}  \sum\limits_{\sigma  \in \mathcal E_{{K^*}}^*} |P_1(v^{(i)}_h)|_{H^2(T_{K^*,\sigma})} + \frac{C_{5}}{2}||\nabla P_1(v_h^{(i)})||_{L^2(K^*)} }{||P_1(\bold  v_h)||_{(H^1(\Omega))^2}}} } } \nonumber \\
&\nonumber \\
&\le & \sqrt{h} \sum\limits_{{K^*} \in {M^*}}  \left(\frac{h}{\text{diam}(K^*)} \right)^{3/2}\frac{(\text{diam}(K^*))^{3/2}}{{(m(K^*)}^{3/4}} \frac{C_{5}}{\sqrt[4]{m(K^*)}} \int\limits_{K^*} |q_h| \text{d}\bold x { { \sum\limits_{i = 1}^2 {\frac{ ||\nabla P_1(v_h^{(i)})||_{L^2(K^*)} }{||P_1(\bold  v_h)||_{(H^1(\Omega))^2}} } } } \nonumber \\
&\nonumber \\
&&\text{because $P_1(v^{(i)}_h)$ is a polynomial of degree $1$, this implies $|P_1(v^{(i)}_h)|_{H^2(T_{K^*,\sigma})} = 0$.} \nonumber \\
&\nonumber \\
& \le & (\sqrt{\zeta_{\mathcal D^{**}}} C_2) ^3\sqrt{h}\left(\sum\limits_{{K^*} \in {M^*}} \sqrt{m(K^*)}\right)^{1/2}
\left(\sum\limits_{{K^*} \in {M^*}} ||q_h||^2_{L^2(K^*)} \sum\limits_{i = 1}^2\frac{||\nabla P_1(v^{(i)}_h) ||^2_{L^2(K^*)}}{||P_1(\bold  v_h)||^2_{(H^1(\Omega))^2}} \right)^{1/2} \nonumber \\
&\nonumber \\
&\le & \underbrace \xi_{=(\sqrt{\zeta_{\mathcal D^{**}}} C_2) ^3}   \sqrt{h}\left(\sum\limits_{{K^*} \in {\mathcal M^*}} \sqrt{m(K^*)}\right)^{1/2} \frac{||\nabla P_1(\bold v_h)||_{(L^2(\Omega))^2}}{||P_1(\bold  v_h)||_{(H^1(\Omega))^2}} ||q_h||_{L^2(\Omega)}.
\end{eqnarray}
From the two inequalities (\ref{stab2}) and (\ref{stab4}), we get 
\begin{eqnarray}
\label{stab5}
&&\mathop {\sup }\limits_{\bold{v}_h \in {\mathcal H_{{\mathcal D}}}\hfill\atop
{\bold v_h} \ne 0\hfill}\int\limits_\Omega  { \frac{\left({\text{div}_{{\mathcal D,\mu }}{\bold v_h}} \right)}{||P_1(\bold v_h)||_{(H^1(\Omega))^2}} {q_h}\,{\text{d}}\bold x } \ge \left[ \eta- \xi\sqrt{h}\left(\sum\limits_{{K^*} \in {\mathcal M^*}} \sqrt{m(K^*)}\right)^{1/2} \right] ||q_h||_{L^2(\Omega)} \nonumber\\
\end{eqnarray}
We assume that there exist $h_0$ be small enough, such that $h \le h_0$
\begin{displaymath}
{h}\left(\sum\limits_{{K^*} \in {\mathcal M^*}} \sqrt{m(K^*)}\right)^{1/2} < \left(\frac{\eta}{\xi}\right )^2.\quad \square
\end{displaymath}
\section{Convergence of the pFECC scheme}
\label{sec5}
In this section we prove that the pair discrete solution $(\bold u_h,p_h) \in \mathcal H_{\mathcal D_h} \times \mathcal L_{\mathcal D}$ tend to the weak solutions $(\bold u, p)$ of the problem (\ref{weakptr2}), as $h \to 0$. \\
We firstly state the theorem \ref{sec5}.1  to prove the convergence of the  velocity.\\
\newline
{\bf Theorem \ref{sec5}.1 (the convergence of the velocity)} 
\textit{ Under hypotheses $\bold f \in (L^2(\Omega))^2$, (\ref{ptr2}) and (\ref{ptr3}), let the positive parameter $\lambda$ be fixed, then  $\bold u_h$ converges to $\bold u$} in $(L^2(\Omega))^2$.\\
\newline
\textit{Proof.}\\
We will prove there exists a sub-sequence of $(\bold u_h)$, such that this sub-sequence converges to $\bold u \in (H^1_0(\Omega))^2$, as $h \to 0$. For this purpose, in Eq.(\ref{disweakptr2a}) and (\ref{disweakptr2b}), we choose $\bold v_h = \bold u_h \in \mathcal H_{\mathcal D}$ and $p_h = q_h \in \mathcal L_{\mathcal D}$. The two equations are rewritten by
\begin{eqnarray}
\label{ptr38}
\int\limits_\Omega  {\mu (\bold x)\nabla_{\mathcal D,\mu}  u^{(1)}_h.\nabla_{\mathcal D,\mu}  u^{(1)}_h\,d\bold x} &+& \int\limits_\Omega  {\mu (\bold x)\nabla_{\mathcal D,\mu}  u^{(1)}_h.\nabla_{\mathcal D,\mu}  u^{(1)}_h\,dx} \nonumber \\
&-& \int\limits_\Omega  {\text{div}_{\mathcal D,\mu} (\bold u_h)\,p_h\,d\bold x}  = \int\limits_\Omega  {\bold f.P(\bold u_h)\,\text{d}\bold x},
\end{eqnarray}
and
\begin{equation}
\label{ptr39}
\int\limits_\Omega  {\text{div}_{\mathcal D,\mu}(\bold u_h)\,p_h\,d\bold x}  = -\lambda h \int\limits_\Omega {p_h}^2 \text{d} \bold x.
\end{equation}
On the left hand side of (\ref{ptr38}), we transform the integral depended the discrete pressure $p_h$ by (\ref{ptr39}), then we get
\begin{eqnarray}
\label{ptr40}
\int\limits_\Omega  {\mu (\bold x)\nabla_{\mathcal D,\mu}  u^{(1)}_h.\nabla_{\mathcal D,\mu}  u^{(1)}_h\,d\bold x} &+& \int\limits_\Omega  {\mu (\bold x)\nabla_{\mathcal D,\mu}  u^{(1)}_h.\nabla_{\mathcal D,\mu}  u^{(1)}_h\,\text{d}\bold x} \nonumber \\
&+& \lambda h \int\limits_\Omega {p_h}^2 \text{d} \bold x = \int\limits_\Omega  {\bold f. (P(u^{(1)})_h,P(u^{(2)}_h))\,\text{d}\bold x}.
\end{eqnarray}
By the condition (\ref{ptr2}) of the  viscosity $\mu$, the left hand side (LHS) of Equation (\ref{ptr40}) is estimated by
\begin{equation}
\label{ptr41}
\text{LHS} \ge \underline{\lambda} \left(||\nabla_{\mathcal D,\mu}u^{(1)}_h||^2_{(L^2(\Omega))^2} +  ||\nabla_{\mathcal D,\mu}u^{(2)}_h||^2_{(L^2(\Omega))^2}\right) + \lambda h_{\mathcal D^{**}} ||p_h||^2_{L^2(\Omega)}.
\end{equation}
And its right hand side (RHS) is bounded by
\begin{eqnarray}
\label{ptr42}
\text{RHS} & \mathop  \le \limits_{{\text{Young}}} & \beta ||\bold f||^2_{(L^2(\Omega))^2} + \frac{{1}}{{\beta}} \left(||P (u^{(1)}_h)||^2_{(L^2(\Omega))} +  ||P (u^{(1)}_h)||^2_{(L^2(\Omega))}\right) \nonumber \\
& \le &  \beta ||\bold f||^2_{(L^2(\Omega))^2} + 
\frac{{C_8}}{{\beta}}
\left(||\nabla_{\mathcal D,\mu}u^{(1)}_h||^2_{(L^2(\Omega))^2} +  ||\nabla_{\mathcal D,\mu}u^{(2)}_h||^2_{(L^2(\Omega))^2}\right)
 \end{eqnarray}
because of the inequality (30) in \cite{PoHa}, where the positive constant $\beta$ is chosen in (\ref{ptr42c}).\\
From the two inequalities (\ref{ptr41}) and (\ref{ptr42}), we have
\begin{equation}
\label{ptr42a}
\left(\underline{\lambda}-\frac{{C_9}}{{\beta}}\right) \left(||\nabla_{\mathcal D,\mu}u^{(1)}_h||^2_{(L^2(\Omega))^2} +  ||\nabla_{\mathcal D,\mu}u^{(2)}_h||^2_{(L^2(\Omega))^2}\right) +  \lambda h||p_h||^2_{L^2(\Omega)} \le \beta ||\bold f||^2_{(L^2(\Omega))^2}
\end{equation}
Using additionally Inequality (21) in \cite{PoHa}, we can estimate Inequality (\ref{ptr42a}) in the discrete $H^1$ norm $||.||_{1,\mathcal D^{**}}$, as follows
\begin{equation}
\label{ptr42b}
\left(\underline{\lambda}-\frac{{C_9}}{{\beta}}\right)C_{10} \left(||u_h^{(1)}||^2_{1,\mathcal D^{**}} +  ||u_h^{(2)}||^2_{1,\mathcal D^{**}}\right) +  \lambda h ||p_h||^2_{L^2(\Omega)} \le \beta ||\bold f||^2_{(L^2(\Omega))^2},
\end{equation}
where $C_9$, $C_{10}$ are not depend on $h$.\\
\newline
Note that the two coefficients $C_9$, $C_{10}$ are generated from the two inequalities $(21)$, $(30)$ in \cite{PoHa}. Besides, the coefficient $\beta$ is chosen by
\begin{equation}
\label{ptr42c}
\beta = \frac{2C_9}{\underline{\lambda}} > 0
\end{equation}
Hence, Inequality (\ref{ptr42b}) follows
\begin{equation}
\label{ptr43}
\left(||u^{(1)}_h||^2_{1,\mathcal D^{**}}+||u^{(2)}_h||^2_{1,\mathcal D^{**}}\right) \le \frac{2\beta}{\underline \lambda}||\bold f||^2_{(L^2(\Omega))^2}
\end{equation}
and
\begin{equation}
\label{ptr43a}
\lambda h||p_h||^2_{L^2(\Omega)} \le \beta||\bold f||^2_{(L^2(\Omega))^2}
\end{equation}
With Inequality $(\ref{ptr43})$, we obtain the existence of a subsequence of $\bold u_h$ and $\bold u \in  (H^1_0(\Omega))^2$ such that this subsequence of $\bold u_h$ converges to $\bold u$ in $(L^2(\Omega))^2$ as $h \to 0$, which is implied from Lemma $5.7$ of \cite{EGH1}.\\
\newline
Let $ \boldsymbol { \phi} \in (C^{\infty}_c(\Omega))^2$ such that $\text{div}(\boldsymbol \phi) = 0$. We suppose that $h$ is small enough, so that, for all $K^{**} \in \mathcal M^{**}$, the  intersection of the two sets $K^{**}$ and $supp \{\boldsymbol \phi\}$ is nonempty, then $\partial K^{**} \cap \partial \Omega = \emptyset$. Additionally, we will give $\bold v_h = \boldsymbol \phi_h  = (\{\boldsymbol \phi(\bold x_{K})\}_{K \in \mathcal M},\{\boldsymbol \phi(\bold x_{K^*})\}_{K^* \in \mathcal M^*}) \in \mathcal H_{\mathcal D}$ in two equations (\ref{disweakptr2a}) and (\ref{disweakptr2b}).
And using the results of Section 5 in \cite{PoHa}, they help us show the convergence of the diffusion operator with the variable viscosity $\mu(\bold x)$:
\begin{eqnarray}
\label{ptr43b}
\mathop {\lim }\limits_{{h} \to 0}  && \left( {\int\limits_\Omega  {\mu (x){\nabla _{\mathcal D,\mu }}u_h^{(1)}.{\nabla _{\mathcal D,\mu }}\phi _h^{(1)}\,{\rm{d}}{\bf{x}}}   +  \int\limits_\Omega  {\mu (x){\nabla _{\mathcal D,\mu }}u_h^{(2)}.{\nabla _{\mathcal D,\mu }}\phi _h^{(2)}\,{\rm{d}}{\bf{x}}} } \right) \nonumber \\
&=& \int\limits_\Omega  {\mu (x)\nabla \bold u:\nabla \boldsymbol \phi \,{\rm{d}}{\bf{x}}} 
\end{eqnarray} 
Furthermore, we apply the Holder inequality and a result $\mathop {\lim }\limits_{h \to 0} {S_i}(\phi ) = 0$ in Corollary of \cite{PoHa} to show  
\begin{equation}
\label{ptr44}
\mathop {\lim }\limits_{{h} \to 0}
\left( {\int\limits_\Omega {\bold f. \boldsymbol \phi_h} \text{d} \bold x}\right) = \int\limits_\Omega {\bold f. \boldsymbol \phi} ~\text{d} \bold x
\end{equation}
Next, we need to prove 
\begin{equation}
\label{ptr45}
\mathop {\lim }\limits_{{h} \to 0}
\left( {\int\limits_\Omega {p_h \text{div}_{\mathcal D, \mu} \boldsymbol \phi_h} \text{d} \bold x}\right) = \mathop {\lim }\limits_{{h} \to 0}
\left[ {\int\limits_\Omega {p_h  \left(\nabla_{\mathcal D, \mu}  \phi^{(1)}_h.\bold e_1 + \nabla_{\mathcal D, \mu}  \phi^{(2)}_h.\bold e_2 \right )}\text{d} \bold x}\right] = 0.
\end{equation}
For this purpose,  we put
\begin{equation}
\label{ptr46}
G^{(i)}_h = \int\limits_\Omega {p_h  (\nabla_{\mathcal D, \mu}  \phi^{(i)}_h.\bold e_i)}\text{d} \bold x = \sum\limits_{{K^*} \in {\mathcal M^*}} {\sum\limits_{{L^{**}} \in {\mathcal M}_{{K^*}}^{**}} {m\left( {{L^{**}}} \right){p_{{K^*}}}\left( {{\nabla _{\mathcal D,\mu }}\phi _h^{(1)}.{e_i}} \right)} },
\end{equation}
\begin{equation}
\label{ptr47}
\overline{G}^{(i)} = \int\limits_\Omega  {{p_h}\left( {\nabla {\phi ^{(i)}}.{{\bf{e}}_i}\,} \right){\rm{d}}{\bf{x}}},
\end{equation}
where their relation is expressed by
\begin{equation}
\label{ptr48}
\mathop {\lim }\limits_{h \to 0} \left( {G_{h}^{(i)} - \overline{G}^{(i)}} \right) = 0 \quad \text{for each}~i=1,2.
\end{equation}
This result is shown by Lemma \ref{sec5}.1, while its condition (\ref{lem2.1})  is satisfied by (\ref{ptr43a}).\\
Besides, we obtain $\overline G^{(1)} + \overline G^{(2)} = 0$ because of $\text{div} \boldsymbol \phi = 0$, which implies (\ref{ptr45}) is proven.\\
\newline
In the last requirement for proving this theorem, we also need to indicate $\text{div}(\bold u) = 0$ a.e in $\Omega$.\\
Let us $\varphi \in C^{\infty}_c(\Omega)$, and the characteristic function $\varphi_h \in \mathcal L_{\mathcal D_h}$ be defined by the value $\varphi(\bold x_{K^*})$ in $K^*$, for all $K^* \in \mathcal M^*$. Remark that the sequence $\varphi_h \to \varphi$ in $L^2(\Omega)$, as $h \to 0$, which is proven by
\begin{equation}
\label{ptr48a}
\int\limits_\Omega |\varphi_h - \varphi(\bold x)|^2 \text{d}\bold x = \sum\limits_{{K^*} \in {\mathcal M^*}} \int\limits_{K^{*}} |\varphi(x_{K^*}) - \varphi(\bold x)|^2 \text{d}\bold x \le C_{11} \sqrt{h}~\text{m}(\Omega) ||\nabla \varphi||_{L^2(\Omega)}.
\end{equation}
To get the inequality (\ref{ptr48a}), we thank to Theorem 9.12 (Morrey) in \cite{Bre}, $\varphi \in C^{\infty}_c(\Omega)$ and $|\bold x - \bold x_{K^*}| \le h$ for all $\bold x \in K^*$ ($\bold x$ must belong to a triangle $T_{K^*} \in \mathcal M^{**}_{K^*}$), $K^* \in \mathcal M^*$.\\
\newline
In Equation (\ref{disweakptr2b}), we choose $q_h = \varphi_h$. This equation is then rewritten by
\begin{equation}
\label{ptr49}
\int\limits_\Omega  {\text{div}_{\mathcal D,\mu}(\bold u_h)\,\varphi_h\,d\bold x}  = -\lambda h \int\limits_\Omega {p_h \varphi_h} \text{d} \bold x,
\end{equation}
whose right and left hand sides are put
\begin{equation}
\label{ptr51}
\overline H_h = \lambda h \int\limits_\Omega {p_h \varphi_h} \text{d} \bold x
\end{equation}
and
\begin{equation}
\label{ptr50}
H_h = \int\limits_\Omega  {\text{div}_{\mathcal D,\mu}(\bold u_h)\,\varphi_h\,d\bold x} = \int\limits_\Omega  {(\nabla_{\mathcal D,\mu} u^{(1)}_h.\bold e_1+\nabla_{\mathcal D,\mu} u^{(2)}_h.\bold e_2)\,\varphi_h\,d\bold x} 
\end{equation}
For the right hand side $\overline H_h$, thanks to (\ref{ptr43a}), it follows
\begin{equation}
\label{ptr52}
|\overline H_h| \mathop  \le \limits_{\scriptstyle\,{\text{Holder}}\hfill} \lambda h ||p_h||_{L^2(\Omega)}||\varphi_h||_{L^2(\Omega)} \le 
 \sqrt{\lambda h} \beta ||f||_{L^2(\Omega)}||\varphi||_{L^2(\Omega)},
\end{equation}
hence, we get 
\begin{equation}
\label{ptr53}
\mathop {\lim }\limits_{{h} \to 0} {\overline H_h} = 0
\end{equation}
For the left hand side $ H_h$, for each direction $x_i$, $i = 1,2$, we use the triangular inequality to evaluate
\begin{equation}
\label{lem3.2}
\left|\int\limits_\Omega  {(\nabla_{\mathcal D,\mu} u^{(i)}_h.\bold e_i)\,\varphi_h\,d\bold x} + \int\limits_\Omega u^{(i)}(\bold x) (\nabla \varphi(\bold x). \bold e_i) \right| \le H_{1,h} + H_{2,h} + H_{3,h},
\end{equation}
where the notations $H_{j,h}$, for $j = 1, 3$, are defined by
\begin{eqnarray*}
H_{1,h} &=& \left| \int\limits_\Omega  \left( \nabla_{\mathcal D,\mu} u_h^{(i)}.\bold e_i\right) (\varphi_h  - \varphi) \text{d} \bold x \right| \\
H_{2,h} &=&  \left| \int\limits_\Omega (\nabla _{\mathcal D,\mu} u_h^{(i)}(\bold x) .\bold e_i) \varphi \text{d} \bold x + \int\limits_\Omega   P(u_h^{(i)}) (\bold x) (\nabla \varphi.\bold e_i) \text{d} \bold x \right| \\
H_{3,h} &=& \left|\int\limits_\Omega u^{(i)}(\bold x) (\nabla \varphi (\bold x).\bold e_i) \text{d} \bold x - \int\limits_\Omega P(u_h^{(i)})(\bold x) (\nabla \varphi (\bold x) .\bold e_i) \text{d} \bold x \right| \\
\end{eqnarray*}
Two coefficients $H_{1,h}$, $H_{3,h}$ are estimated by the Holder inequality
\begin{equation}
\label{ptr54}
H_{1,h} \le ||\nabla_{\mathcal D,\mu} u^{(i)}_h.\bold e_i||_{L^2(\Omega)} ||\varphi_h - \varphi||_{L^2(\Omega)} \le \frac{2\beta}{\underline \lambda} ||\varphi_h - \varphi||_{L^2(\Omega)}
\end{equation}
and 
\begin{equation}
\label{ptr55}
H_{3,h} \le ||\nabla \varphi.\bold e_i||_{L^2(\Omega)} ||u^{(i)}_h - u^{(i)}||_{L^2(\Omega)}.
\end{equation}
From the above results, we claim that $H_{1,h}$, $H_{3,h}$ tend to $0$, while the sequence $\varphi_h \to \varphi$ and $u_h^{(i)} \to u^{(i)}$ in $L^2(\Omega)$, as $ h \to 0$. \\
\newline
Now, we have to prove $H_{2,h} \to 0$, as $h \to 0$. 
Before, let us introduce the two following sets
\begin{displaymath}
\mathcal M^{**}_{\mu} = \left\{  K \in \mathcal M^{**} | \mu ~\text{is not continuous on}~ K \right\},
\end{displaymath}
if $K \in \mathcal M^{**} \ \mathcal M^{**}_{\text{const}}$, then we assume that $\mu(\bold x)= \mu_1$ on a triangle $K_1 =(\bold x_K, \bold x_{K^*}, \bold x_\sigma)$, that $\mu(\bold x)= \mu_2$ on a triangle $K_2 = (\bold x_L, \bold x_{K^*}, \bold x_\sigma)$. And
\begin{displaymath}
\mathcal M^{**}_{\text{const}} = \left\{ K \in \mathcal M^{**} | \Lambda ~\text{is constant on}~ K \right \} 
\end{displaymath}
We rewrite 
\begin{displaymath}
H_{2,h} = {\hat H}_{2,h} + \overline{\overline{H}}_{2,h}, 
\end{displaymath}
where 
\begin{eqnarray*}
{\hat H}_{2,h} &=&\int\limits_\Omega (\nabla _{\mathcal D,\mu} u_h^{(i)}(\bold x) .\bold e_i) \varphi \text{d} \bold x +  \int\limits_\Omega   P_1(u_h^{(i)}) (\bold x) (\nabla \varphi.\bold e_i) \text{d} \bold x \\
\overline{\overline{H}}_{2,h} &=& \int\limits_\Omega   [P(u_h^{(i)})-P_1(u_h^{(i)})] (\bold x) (\nabla \varphi.\bold e_i) \text{d} \bold x.
\end{eqnarray*}
We also introduce some notations, as follows: $\varphi_K$ the average value of $\varphi$ if $K \in \mathcal M^{**} \cap \mathcal M^{**}_\mu$, $\varphi_{K_1}$ (resp.$\varphi_{K_2}$) the average value of $\varphi$ on $K_1$ (resp. $K_2$) if $K^{**} \in \mathcal M^{**}_{\mu}$, $\varphi_{M,N}$ with $(M,N) \in S_K$, the average value of $\varphi$ on $\overrightarrow{\tau}_{M,N,(M,N)\in S_K}$. We express $\hat H_{2,h}$ by the sum $ \hat H^{(1)}_{2,h} + \hat H^{(2)}_{2,h} + \hat H^{(3)}_{2,h}$ defined by
\begin{eqnarray*}
\hat H^{(1)}_{2,h} &=& \sum\limits_{K \in \mathcal M_{const}^{**}} {|K|\nabla_K u^{(i)}_h. (\varphi_K \bold e_i)} \nonumber \\
&+& \sum\limits_{K \in \mathcal M^{**} \backslash \{\mathcal M_{const}^{**} \cup \mathcal M^{**}_{\Lambda} \}} \left[ {|K_1|\nabla_{K_1} u^{(i)}_h. (\varphi_{K_1} \bold e_i) + |K_2|\nabla_{K_2} u^{(i)}_h. (\varphi_{K_2} \bold e_i)} \right]  \\
&+& \sum\limits_{K \in \mathcal M_{\Lambda}^{**}} \left[ {|K_1|\nabla_{K_1} u^{(i)}_h. (\varphi_{K_1} \bold e_i) + |K_2|\nabla_{K_2} u^{(i)}_h. (\varphi_{K_2} \bold e_i)} \right]  
\end{eqnarray*}
\begin{eqnarray*}
\hat H^{(2)}_{2,h} &=& \sum\limits_{K \in \mathcal M^{**}} 
\left[ \begin{array}{l}
\vspace{0.4cm}
(u^{(i)}_K - u^{(i)}_{K^*}) \overrightarrow{\tau}_{K,K^*}.(\varphi_{K,K^*} \bold e_i) + (u^{(i)}_{K^*} - u^{(i)}_{L}) \overrightarrow{\tau}_{K^*,L} \\
+ (u^{(i)}_L - u^{(i)}_{K}) \overrightarrow{\tau}_{L,K}.(\varphi_{K,L} \bold e_i).(\varphi_{K,L} \bold e_i)
\end{array} \right] \\
&\\
&\\
\hat H^{(3)}_{2,h} &=&  \sum\limits_{K \in \mathcal M^{**}} 
\left[ \begin{array}{l}
\vspace{0.4cm}
\int\limits_{{A_{{i_K}}}} {\left[ {{\nabla _{{P_{1,K}}}}{u^{(i)}}.(x - {x_K})} \right]\left( {\nabla \varphi .{{\bf{e}}_i}} \right)} {\rm{d}}{\bf{x}}\\
\vspace{0.4cm}
\int\limits_{{A_{{i_L}}}} {\left[ {{\nabla _{{P_{1,L}}}}{u^{(i)}}.(x - {x_L})} \right]\left( {\nabla \varphi .{{\bf{e}}_i}} \right)} {\rm{d}}{\bf{x}}\\
\int\limits_{{A_{{i_{K^*}}}}} {\left[ {{\nabla _{{P_{1,K^*}}}}{u^{(i)}}.(x - {x_{K^*}})} \right]\left( {\nabla \varphi .{{\bf{e}}_i}} \right)} {\rm{d}}{\bf{x}}
\end{array} \right]
\end{eqnarray*}
Using the computational results of $T_1$, $T_2$, $T_3$ represented in the pages $27$, $28$ of \cite{PoHa}, we obtain
\begin{equation}
\label{ptr56}
|\hat H^{(1)}_{2,h} + \hat H^{(2)}_{2,h}| \le C_{12} ||u||_{1,\mathcal D^{**}} (h_{\mathcal D^{**} + \epsilon_2(h)})
\end{equation}
with $\mathop {\lim }\limits_{{h} \to 0} {\varepsilon _2}({h}) = 0$,
\begin{equation}
\label{ptr57}
\hat H^{(3)}_{2,h}\le h ||\nabla_{P_1} u^{(i)}_h||_{(L^2(\Omega))^2} ||\nabla \varphi.\bold e_i||_{L^2(\Omega)}.
\end{equation}
With $\overline{\overline{H}}_{2,h}$, we also use the results (28) of \cite{PoHa} and (\ref{ptr42a}) to get
\begin{equation}
\label{ptr58}
\overline{\overline{H}}_{2,h} \le
 ||P(u^{(i)}_h) - P_1(u^{(i)}_h)||_{L^2(\Omega)} ||\nabla \varphi.\bold e_i||_{L^2(\Omega)} \le C_{13} h ||\bold f||_{(L^2(\Omega))^2} ||\nabla \varphi.\bold e_i||_{L^2(\Omega)}
\end{equation}
From (\ref{ptr56})-(\ref{ptr58}), we obtain
\begin{equation}
\label{ptr59}
H_{2,h} \to 0, \quad \text{as}~h \to 0,
\end{equation}
together the convergence $H_{1,h}$, $H_{3,h}$ to $0$, we conclude that
\begin{equation}
\label{ptr60}
\int\limits_\Omega  {(\nabla_{\mathcal D,\mu} u^{(1)}_h.\bold e_1+\nabla_{\mathcal D,\mu} u^{(2)}_h.\bold e_2)\,\varphi_h\,d\bold x} \to - \int\limits_\Omega \bold u(\bold x).\nabla \varphi(\bold x) = \int\limits_\Omega \text{div} (\bold u) (\bold x) \varphi(\bold x),  
\end{equation}
as $h \to 0$.\\
Therefore, the results (\ref{ptr49}), (\ref{ptr53}) and (\ref{ptr60}) imply 
\begin{equation}
\label{ptr60a}
\int\limits_\Omega \text{div} \bold u(\bold x) \varphi(\bold x) = 0, \quad \text{for all}~\varphi \in C^{\infty}(\Omega). 
\end{equation}
From the results (\ref{ptr43b}), (\ref{ptr44}), (\ref{ptr45}) and (\ref{ptr60a}), we proved that the approximate solution $\bold u_h$ converges to the weak analysis solution $\bold u$.  $\hspace{7cm} \square$\\
\newline
\newline
{\bf Theorem \ref{sec5}.2 (the convergence of the pressure)} 
\textit{ Under hypotheses $\bold f \in (L^2(\Omega))^2$, (\ref{ptr2}) and (\ref{ptr3}), let the positive parameter $\lambda$ be fixed, then, the approximate pressure $p_h$ converges to $p$ in $L^2(\Omega)$}.\\
\newline
\textit{Proof}\\
In Equation (\ref{disweakptr2a}), let any $\boldsymbol {\overline \phi} \in (C^2(\Omega) \cap H_0^1(\Omega))^2 $, we choose
\begin{displaymath}
\bold v_h = \boldsymbol {\overline \phi}_h = \left(\{\boldsymbol {\overline \phi}(\bold x_K)\}_{K\in \mathcal M},\boldsymbol {\overline \phi}(\bold x_{K^*})\}_{K^*\in \mathcal M^*} \right),
\end{displaymath}
so these equations are rewritten by
\begin{equation}
\label{ptr61}
\int\limits_{{\Omega}} {\left( {\mu (\bold x){\nabla _{\mathcal D,\mu }}{(\bold u_h)}} \right):{\nabla _{\mathcal D,\mu }} \boldsymbol {\overline \phi}_h\text{d} \bold x}  - \int\limits_{{\Omega}} {\left( { \text{div}_{\mathcal D,\mu} \boldsymbol {\overline \phi}_h } \right)p_h~\text{d} \bold x}  =  \int\limits_{\Omega} {{\bold f}.{ P(\boldsymbol {\overline \phi}_h)}\,\text{d}{\bold x} }
\end{equation}
and in Equation (\ref{weakptr2}), $\bold v = \bold P_1(\boldsymbol {\overline \phi}_h)$
\begin{equation}
\label{ptr62}
\int\limits_\Omega  {\mu (\bold x)\nabla \bold u:\nabla  \bold P_1(\boldsymbol {\overline \phi}_h)\,d\bold x}  - \int\limits_\Omega  {\text{div} (\bold P_1(\boldsymbol {\overline \phi}_h))\,p\,dx}  = \int\limits_\Omega  {\bold f. P_1(\boldsymbol {\overline \phi}_h)\,d\bold x},
\end{equation}
Equation (\ref{ptr61}) is subtracted to (\ref{ptr21}) equals 
\begin{eqnarray*}
&&\int\limits_{{\Omega}} {\left( {\mu (\bold x){\nabla _{\mathcal D,\mu }}{\bold u_h}} \right):{\nabla _{\mathcal D,\mu }} \boldsymbol {\overline \phi}_h\text{d} \bold x} - \int\limits_\Omega  {\mu (\bold x)\nabla \bold u:\nabla   P_1(\boldsymbol {\overline \phi}_h)\,d\bold x} \\
&&- \int\limits_{{\Omega}} {\left( { \text{div}_{\mathcal D,\mu} \boldsymbol {\overline \phi}_h } \right)p_h~\text{d} \bold x} + \int\limits_\Omega  {\text{div} ( P_1(\boldsymbol {\overline \phi}_h))\,p\,dx}  = \int\limits_\Omega  {\bold f.( P(\boldsymbol {\overline \phi}_h)-  P_1(\boldsymbol {\overline \phi}_h)) \,\text{d}\bold x}.\\
\end{eqnarray*}
This equation corresponds to the following equation
\begin{equation}
\label{ptr63}
O_{1,h} + O_{2,h} + O_{3,h} + O_{4,h} + O_{5,h} = O_{6,h},
\end{equation}
where $P_1(\boldsymbol {\overline \phi}_h)$ is different from $0$, $O_{i,h}$, $i = \overline{1,6}$, are defined, as follows:
\begin{eqnarray*}
O_{1,h} &=& \frac{1}{|| P_1(\boldsymbol {\overline \phi}_h)||_{(H^1(\Omega))^2}}\int\limits_{{\Omega}} {\left( {\mu (\bold x){\nabla _{\mathcal D,\mu }}{\bold u_h}} \right):({\nabla _{\mathcal D,\mu }} \boldsymbol {\overline \phi}_h} - \nabla \boldsymbol {\overline \phi}) \text{d} \bold x, \\
O_{2,h} &=& \frac{1}{|| P_1(\boldsymbol {\overline \phi}_h)||_{(H^1(\Omega))^2}}\int\limits_\Omega  {(\mu (\bold x)\nabla \bold u) :( \nabla \boldsymbol {\overline \phi} - \nabla  \bold P_1(\boldsymbol {\overline \phi}_h) )\text{d}\bold x}, \\
O_{3,h} &=& \frac{1}{|| P_1(\boldsymbol {\overline \phi}_h)||_{(H^1(\Omega))^2}}\int\limits_\Omega  {(\mu (\bold x)\nabla\boldsymbol {\overline \phi}) :( \nabla_{\mathcal D,\mu} \bold u_h - \nabla  \bold u )\text{d}\bold x},\\
O_{4,h} &=& \frac{1}{|| P_1(\boldsymbol {\overline \phi}_h)||_{(H^1(\Omega))^2}}\int\limits_{{\Omega}} {\left( { \text{div}  P_1(\boldsymbol {\overline \phi}_h) } - \text{div}_{\mathcal D,\mu}(\Pi_h \boldsymbol {\overline \phi}) \right) p~\text{d} \bold x},\\
O_{5,h} &=&\frac{1}{|| P_1(\boldsymbol {\overline \phi}_h)||_{(H^1(\Omega))^2}}\int\limits_\Omega  {\bold f.[ P_1(\boldsymbol {\overline \phi}_h)-  P(\boldsymbol {\overline \phi}_h)] \,\text{d}\bold x},\\
O_{6,h} &=& \frac{1}{|| P_1(\boldsymbol {\overline \phi}_h)||_{(H^1(\Omega))^2}}\int\limits_{{\Omega}} {\left( { \text{div}_{\mathcal D,\mu} \boldsymbol {\overline \phi}_h } \right) (p_h-p)~\text{d} \bold x}.
\end{eqnarray*}
Using the Korn inequality, we get
\begin{eqnarray*}
|O_{1,h}| & \le & \overline \lambda \frac{||{{\nabla _{\mathcal D,\mu }}{\bold u_h}}||_{(L^2(\Omega))^2}}{|| P_1(\boldsymbol {\overline \phi}_h)||_{(H^1(\Omega))^2}} ~ || {\nabla _{\mathcal D,\mu }} \boldsymbol {\overline \phi}_h - \nabla \boldsymbol {\overline \phi} ||_{(L^2(\Omega))^2},\\
&\\
|O_{2,h}| &\le & \overline \lambda \frac{||\nabla \bold u ||_{(L^2(\Omega))^2}}{|| P_1(\boldsymbol {\overline \phi}_h)||_{(H^1(\Omega))^2}}~ || \nabla \boldsymbol {\overline \phi} - \nabla   P_1(\boldsymbol {\overline \phi}_h ||_{(L^2(\Omega))^2},\\
&\\
|O_{3,h}| & \le & \overline \lambda \frac{||{{\nabla _{\mathcal D,\mu }}{\bold u_h}}||_{(L^2(\Omega))^2}}{|| P_1(\boldsymbol {\overline \phi}_h)||_{(H^1(\Omega))^2}} ||\nabla \boldsymbol {\overline \phi} - \nabla_{\mathcal D,\mu} \Pi_h (\boldsymbol {\overline \phi})||_{(L^2(\Omega))^2} \\
&+& \frac{1}{|| P_1(\boldsymbol {\overline \phi}_h)||_{(H^1(\Omega))^2}}\left| \int\limits_\Omega  {(\mu (\bold x) \nabla_{\mathcal D,\mu} (\boldsymbol {\overline \phi}_h) : \nabla_{\mathcal D,\mu} \bold u_h\text{d}\bold x}- \int\limits_\Omega  {(\mu (\bold x)\nabla\boldsymbol {\overline \phi}) : \nabla  \bold u \text{d}\bold x} \right |, \\
&\\
|O_{4,h}| & \le & \frac{||p||_{(L^2(\Omega))^2}} {|| P_1(\boldsymbol {\overline \phi}_h)||_{(H^1(\Omega))^2}} \left( || \text{div} P_1(\boldsymbol {\overline \phi}_h) - \text{div} \boldsymbol {\overline \phi} ||_{(L^2(\Omega))^2} +  || \text{div} \boldsymbol{\overline \phi} - \text{div}_{\mathcal D,\mu}(\Pi_h \boldsymbol {\overline \phi})||_{(L^2(\Omega))^2} \right)\\
&\\
|O_{5,h}| & \le & \frac{||\bold f||_{(L^2(\Omega))^2}}{|| P_1(\boldsymbol {\overline \phi}_h)||_{(H^1(\Omega))^2}}~|| P_1( \boldsymbol {\overline \phi}_h)-  P(\boldsymbol {\overline \phi}_h)||_{(L^2(\Omega))^2}.\\
\end{eqnarray*}
When $h_{\mathcal D^{**}}$ tends to $0$,  we have
\begin{itemize}
\item $|| {\nabla _{\mathcal D,\mu }} \boldsymbol {\overline \phi}_h - \nabla \boldsymbol {\overline \phi} ||_{(L^2(\Omega))^2} \to  0 \quad \text{because of Lemma $4.3$ in \cite{EGH1}},$ 
\item $|| \nabla \boldsymbol {\overline \phi} - \nabla   P_1(\boldsymbol {\overline \phi}_h) ||_{(L^2(\Omega))^2}  \to  0, ~ || \text{div}  P_1(\boldsymbol {\overline \phi}_h) - \text{div} \boldsymbol {\overline \phi} ||_{(L^2(\Omega))^2} \to 0 $\\
\newline
because of Theorems $3.4.3$, $3.4.4$ in \cite{QuVa}.
\item $|| \text{div} \boldsymbol{\overline \phi} - \text{div}_{\mathcal D,\mu}(\boldsymbol {\overline \phi}_h)||_{(L^2(\Omega))^2} \to 0,$ $||P_1( \boldsymbol {\overline \phi}_h)-  P(\boldsymbol {\overline \phi}_h)||_{(L^2(\Omega))^2} \to 0$\\
\newline
and $|| P_1( \boldsymbol {\overline \phi}_h)||_{(L^2(\Omega))^2} \to ||\boldsymbol {\overline \phi}||_{(L^2(\Omega))^2}$ because of Proposition $5.3$, Corollary $5.4$ in \cite{PoHa}
\item $\left | \int\limits_\Omega  {(\mu (\bold x) \nabla_{\mathcal D,\mu} ( \boldsymbol {\overline \phi}_h) : \nabla_{\mathcal D,\mu} \bold u_h\text{d}\bold x}- \int\limits_\Omega  {(\mu (\bold x)\nabla\boldsymbol {\overline \phi}) : \nabla  \bold u \text{d}\bold x} \right | \to 0$\\
\newline
because of the above Lemma $4.1$
\end{itemize}
Moreover, $||\nabla_{\mathcal D, \mu} \bold u_h ||_{(L^2(\Omega))^2}$ is upper bounded by the positive constant independent $h_{\mathcal D^{**}}$, which is implied from (\ref{ptr42a})-(\ref{ptr42c}). \\
\newline
In order to $O_{6,h}$, we rewrite Lemma $4.2$, as follows: let any $\epsilon >0$, there then exists $\bold v_h \in \mathcal H_{\mathcal D}$, such that 
\begin{displaymath}
\frac{{\int\limits_\Omega  {\left( {\text{div}_{{\mathcal D,\mu }}{\bold v_h}} \right){q_h}\,{\text{d}}\bold x } }}{{||P_1(\bold  v_h)||_{(H^1(\Omega))^2}}} + \epsilon \ge C_{14}||{q_h}|{|_{{L^2}(\Omega )}}
\end{displaymath}
With  $\bold v_h \in  \mathcal H_{\mathcal D}$ satisfied the above inequality, we apply Theorem $1$ and Corollary in \cite{Ze} and Theorem $2$ in \cite{Ze1} to generate $\boldsymbol \varphi $ a triangular $C^2$-element of degree $9$ such that $\boldsymbol \varphi \in (C^2(\Omega) \cap H_0^1(\Omega))^2$ and  $\bold v_h = \boldsymbol \varphi_h = \left(\{\boldsymbol \varphi(\bold x_K)\}_{K\in \mathcal M},\boldsymbol \varphi(\bold x_{K^*})\}_{K^*\in \mathcal M^*} \right)$. This implies that
\begin{displaymath}
\int\limits_\Omega  { \frac{{\text{div}_{{\mathcal D,\mu }}{\bold v_h}} }{||P_1( \bold v_h)||_{(H^1(\Omega))^2}} {q_h}\,{\text{d}}\bold x } = \int\limits_\Omega  {\frac{{{\text{div}_{\mathcal D,\mu }}\boldsymbol\varphi_h}}{{||{P_1}(\boldsymbol \varphi_h)|{|_{{{\left( {{H^1}(\Omega )} \right)}^2}}}}}\,{q_h}{\rm{d}}{\bf{x}}}.
\end{displaymath}
Hence,
\begin{displaymath}
\int\limits_\Omega  {\frac{{{\text{div}_{\mathcal D,\mu }}\boldsymbol \varphi_h }}{{||\boldsymbol\varphi_h |{|_{{{\left( {{H^1}(\Omega )} \right)}^2}}}}}\,{q_h}{\rm{d}}{\bf{x}}} + \epsilon \ge C_{14}||{q_h}|{|_{{L^2}(\Omega )}},
\end{displaymath}
this corresponds to 
\begin{equation}
\label{ptr64}
\mathop {\sup }\limits_{\boldsymbol \varphi \in {(C^2(\Omega) \cap H^1_0(\Omega))^2}\hfill\atop
\quad \quad  \boldsymbol \varphi \ne 0\hfill} \int\limits_\Omega  {\frac{{{\text{div}_{D,\mu }}\boldsymbol \varphi_h}}{{||{P_1}(\boldsymbol \varphi_h )|{|_{{{\left( {{H^1}(\Omega )} \right)}^2}}}}}\,{q_h}{\rm{d}}{\bf{x}}} \ge C_{14}||{q_h}|{|_{{L^2}(\Omega )}}.
\end{equation}
From Equations (\ref{ptr63}), (\ref{ptr64}), the above estimations of $|O_{i,h}|$ with $i = \overline{1,5}$ and Inequality (\ref{stab1}), we obtain\\
\begin{equation}
\label{ptr65}
||p_h - p||_{L^2(\Omega)} \to 0, \quad \text{as}~h \to 0.\quad\square
\end{equation}
\newline
\textbf{Lemma \ref{sec5}.1} \textit{For a given sequence $p_h \in \mathcal L_{\mathcal D_h}$, it satisfies 
\begin{equation}
\label{lem2.1}
||p_h||_{L^2(\Omega)} \le \frac{C_{15}}{\sqrt{h}},
\end{equation}
where the constant $C_{15}$ is positive,  then 
\begin{equation}
\label{lem2.2}
\mathop {\lim }\limits_{{h} \to 0} \left(\sum\limits_{{K^*} \in {\mathcal M^*}} {\sum\limits_{{L^{**}} \in {\mathcal M}_{{K^*}}^{**}} {m\left( {{L^{**}}} \right){p_{{K^*},h}}\left( {{\nabla _{\mathcal D,\mu }}{\phi}_h^{(i)}.{e_i}} \right)} } - \int\limits_\Omega  {{p_h}\left( {\nabla {{ \phi} ^{(i)}}.{{\bf{e}}_i}\,} \right){\rm{d}}{\bf{x}}}\right) = 0
\end{equation}
for each $i = 1, 2$, where the above vector function $\boldsymbol \phi = (\phi^{(1)},\phi^{(2)})$ satisfies the same conditions as $\boldsymbol \phi$ in Theorem \ref{sec5}.1, and the integral $\int\limits_\Omega  {{p_h}\left( {\nabla_{\mathcal D, \mu} {{ \phi} ^{(i)}}.{{\bf{e}}_i}\,} \right){\rm{d}}{\bf{x}}}$ is written by the formula
$\sum\limits_{{K^*} \in {\mathcal M^*}} {\sum\limits_{{L^{**}} \in {\mathcal M}_{{K^*}}^{**}} {m\left( {{L^{**}}} \right){p_{{K^*},h}}\left( {{\nabla _{\mathcal D,\mu }}{ \phi} _h^{(i)}.{e_i}} \right)} }$.}\\
\newline
Note that the approximate pressure $p_h$ satisfies the condition (\ref{lem2.1}), because of (\ref{ptr42b}).\\
\newline
\textit{Proof.}\\
We will define the two notations $G^{(i)}_{1,h}$ and $G^{(i)}_{2,h}$, as follows:
\begin{displaymath}
G^{(i)}_{1,h} = 
\sum\limits_{K^* \in {\mathcal M}^*} \sum\limits_{{L^{**}} \in {\mathcal M}_{{K^*}}^{**}} {m\left( {{L^{**}}} \right) {p_{{K^*},h}}\left( {{\nabla _{\mathcal D,\mu }}{ \phi} _h^{(1)}.{e_i}} \right)} ~ \text{and} ~ G^{(i)}_{2,h} = \int\limits_\Omega  {{p_h}\left( {\nabla {{ \phi} ^{(i)}}.{{\bf{e}}_i}\,} \right){\rm{d}}{\bf{x}}}.
\end{displaymath}
We then have their computations
\begin{eqnarray}
\label{lem2.3}
&&|G^{(i)}_{1,h}-G^{(i)}_{2,h}| = \left| \sum\limits_{{K^*} \in {\mathcal M^*}} {{p_{{K^*},h}}\left[ {\sum\limits_{{L^{**}} \in \mathcal M_{{K^*}}^{**}} {m\left( {{L^{**}}} \right)\left( {{\nabla _{\mathcal D,\mu }}{ \phi} _h^{(i)}.{\bold e_i}} \right) - \int\limits_{{K^*}} {\left( {\nabla {{ \phi} ^{(i)}}.{{\bold{e}}_i}\,} \right){\rm{d}}{\bf{x}}} } } \right]}\right| \nonumber \\ 
&\le & \sum\limits_{{K^*} \in {\mathcal M^*}} {\sqrt{\text{m}_{K^*}}\left| {{p_{{K^*},h}}} \right|\frac{1}{\sqrt{\text{m}_{K^*}}}  \left| {\sum\limits_{{L^{**}} \in \mathcal M_{{K^*}}^{**}} {m\left( {{L^{**}}} \right)\left( {{\nabla _{\mathcal D,\mu }}{ \phi} _h^{(i)}.{\bold e_i}} \right) - \int\limits_{\partial {K^*}} {{{ \phi} ^{(i)}}\left( {\bar n({\bf{x}}).{{\bold{e}}_i}\,} \right){\rm{d}}{\bf{x}}} } } \right|}\nonumber\\
\end{eqnarray}
Applying Lemma \ref{sec4}.1 to the term $G^{(i)}_{K^*}$ in the right hand side of the inequality (\ref{lem2.3})
\begin{eqnarray}
\label{lem2.4}
|G^{(i)}_{K^*}| &=& \left| {\sum\limits_{{L^{**}} \in \mathcal M_{{K^*}}^{**}} {m\left( {{L^{**}}} \right)\left( {{\nabla _{\mathcal D,\mu }}{ \phi} _h^{(i)}.{\bold e_i}} \right) - \int\limits_{\partial {K^*}} {{{ \phi} ^{(i)}}\left( {\bar n({\bf{x}}).{{\bold{e}}_i}\,} \right){\rm{d}}{\bf{x}}} } } \right| \nonumber \\
& \le & {C_{4}}h|{{{ \phi}^{(i)}}}{|_{{H^2}({K^*})}}
\end{eqnarray}
and using the following estimation for any $K^* \in \mathcal M^*$
\begin{equation}
\label{lem2.5}
\frac{h}{\sqrt{{m}(K^*)}} \le \frac{h}{\sqrt{{m}(L^{**})}} \le \frac{h}{\text{diam}(L^{**})\sqrt{C_{3}}} \le \frac{\zeta_{\mathcal D^{**}}}{\sqrt{C_{3}}}
\end{equation}
with any $L^{**}\in \mathcal M_{K^*}^{**}$, they lead to 
\begin{eqnarray}
\label{lem2.6}
|G^{(i)}_{1,h}-G^{(i)}_{2,h}| &\le& \sum\limits_{{K^*} \in {\mathcal M^*}} {\sqrt{{m}({K^*})}\left| {{p_{{K^*}}}} \right|\frac{h^2}{\sqrt{{m}({K^*})}} {C_{4}}|{{{ \phi}^{(i)}}}{|_{{H^2}({K^*})}} } \nonumber \\
&\mathop  \le \limits_{\scriptstyle\,{\rm{Cauchy}}\hfill\atop
\scriptstyle{\rm{Schwarz}}\hfill} & ||p_h||_{L^2(\Omega)}. \frac{h\zeta_{\mathcal D^{**}}}{{\sqrt{C_{3}}}}  C_{4} |{{{ \phi}^{(i)}}}{|_{{H^2}({\Omega})}} \mathop  \le \limits_{\,({\ref{lem2.1}})} {C_{16}}{h}|{{ \phi}^{(i)}}|_{{{H^2}(\Omega )}}. \nonumber \\
\end{eqnarray}
Therefore, 
\begin{displaymath}
\mathop {\lim }\limits_{{h} \to 0} \left( {G_{1,h}^{(i)} - G_{2,h}^{(i)}} \right) = 0. 
\end{displaymath}
$\hspace{15cm} \square$

\end{document}